\documentclass[12pt]{article}
\usepackage{amsfonts}
\def\Bbb#1{\mathbb#1}
\newtheorem{theorem}{Theorem}

\newtheorem{claim}[theorem]{Claim}
\newtheorem{proposition}[theorem]{Proposition}
\newtheorem{lemma}[theorem]{Lemma}
\newtheorem{corollary}[theorem]{Corollary}

\newtheorem{examples}[theorem]{Examples}

\newcommand{\R}{\Bbb{R}}
\newcommand{\Q}{\Bbb{Q}}
\newcommand{\Sp}{\Bbb{S}}
\newcommand{\Sf}{\Bbb{S}}

\newcommand{\Les}{\Bbb{L}}

\newcommand{\Hy}{\Bbb{H}}

\newcommand{\spa}{\mbox{span}}
\newcommand{\hess}{\mbox{Hess\,}}

\newcommand{\po}{{\hspace*{-1ex}}{\bf .  }}
\newcommand{\ii}{isometric immersion }

\newcommand{\sff}{second fundamental form }

\newcommand\fall{\;\;\;\mbox{for all}\;\;}

\def\P{{\cal P}}
\def\Nu{{\cal V}}
\def\Ral{{\cal R}}
\def\Sal{{\cal S}}
\def\Fes{{\cal F}}

\newcommand{\E}{{\cal E}}
\newcommand{\U}{{\cal U}}
\newcommand{\N}{{\cal N}}
\newcommand{\Np}{{\cal N\,^{^\|}}}

\def\<{\langle}

\def\>{\rangle}
\def\a{\alpha}
\def\va{\varphi}

\def\o{\omega}

\def\id{I}
\def\e{\epsilon}
\def\d{\partial}
\def\bea{\begin{eqnarray*} }
\def\eea{\end{eqnarray*} }
\def\be{\begin{equation} }
\def\ee{\end{equation} }

\def\nap{\nabla^\perp}

\def\nab{\tilde\nabla}
\def\proof{\noindent{\it Proof: }} 
\def\qed{\ifhmode\unskip\nobreak\fi\ifmmode\ifinner\else
\hskip5 pt \fi\fi\hbox{\hskip5 pt \vrule width4 pt  
height6 pt  depth1.5 pt \hskip 1pt }}

\newcommand{\Lal}{{\cal L}}
\begin{document}

\title{Reducibility of Dupin submanifolds}
\author {Marcos Dajczer, Luis A. Florit and Ruy Tojeiro}
\date{}
\maketitle

\begin{abstract}
We introduce the notion of weak reduciblity for Dupin submanifolds with arbitrary codimension. We give a
complete characterization of all weakly reducible Dupin submanifolds, as a consequence of a general result on a 
broader class of Euclidean submanifolds. As a main application, we derive an explicit recursive procedure to generate all holonomic Dupin submanifolds in terms of solutions of 
completely integrable systems of linear partial differential equations of first order. We obtain several additional results
on Dupin submanifolds.
\end{abstract}

  A hypersurface $f\colon\,M^n\to\Q_c^{n+1}$ of a simply 
connected space form of sectional curvature $c$ is called 
{\em proper Dupin\/} if the number of principal curvatures 
is constant and each one of them is constant along the 
corresponding eigenbundle. These conditions are 
\mbox{invariant} under conformal transformations of the 
ambient space, which makes the theory of Dupin hypersurfaces
essentially the same whether it is considered in Euclidean 
space $\R^{n+1}$, in the sphere $\Sf^{n+1}$ or in 
hyperbolic space $\Hy^{n+1}$.
More generally, the class of proper Dupin hypersurfaces 
in Euclidean 
space $\R^{n+1}$ is invariant under  the Lie sphere group  
generated by the subgroup of conformal (Moebius) 
transformations of $\R^{n+1}$ together with the 
\mbox{$1$-parameter} subgroup of parallel translations 
that transform a hypersurface to its parallel one 
at a fixed distance in the normal direction. 
Two hypersurfaces that differ by a Lie transformation 
are said to be Lie equivalent.

An important class of proper Dupin hypersurfaces in
$\R^{n+1}$ is that of
stereographic projections of  isoparametric hypersurfaces
in $\Sf^{n+1}$.  The latter are abundant and have not yet 
been completely classified, although several
interesting results are known including strong restrictions 
on the numbers of distinct principal curvatures
and their multiplicities; see \cite{th}
for a nice recent survey that also discusses Dupin 
hypersurfaces.

It was observed by Pinkall \cite{pi} that further local 
examples of proper Dupin hypersurfaces in $\R^{n+1}$ having 
any given number of principal curvatures with arbitrarily prescribed multiplicities can be constructed by means of 
one of the following procedures, the last two of which yielding submanifolds that are Lie equivalent. Start with a proper Dupin hypersurface $L^{n-s}$ in $\R^{n-s+1}$, the latter regarded 
as a linear subspace $\R^{n-s+1}\times \{0\}$ of 
$\R^{n+1}$, 
and let $M^n$ be defined as
\begin{itemize}
\item[(i)] $M^n$ is the cylinder $L^{n-s}\times \R^s$; 
\item[(ii)] $M^n$ is 
obtained by rotating  $L^{n-s}$ around an axis $\R^{n-s}\subset\R^{n-s+1}$;
\item[(iii)] $M^n$ is the cylinder $CV^{n-s}\times \R^{s-1}$, where $CV^{n-s}$ is the cone over the inverse image $V^{n-s}\subset\Sf^{n-s+1}\subset\R^{n-s+2}$ of $L^{n-s}$ 
by the stereographic projection;
\item[(iii\hspace{0,1ex}$'$)] $M^n$ is a tube around $L^{n-s}$.
\end{itemize}

  These constructions introduce a new principal curvature $\lambda$ of multiplicity~$s$ that is easily seen to be 
constant along its eigenbundle $E_\lambda$. The other 
principal curvatures of $M^n$ are determined from those of $L^{n-s}$, and they are constant along the corresponding eigenbundles because $L^{n-s}$ is Dupin. Moreover, the {\em conullity\/} distribution $E_\lambda^\perp$ of $\lambda$, 
that is, the orthogonal distribution to  $E_\lambda$ in the tangent bundle $TM$, is always integrable. In fact, in the 
first three constructions $E_\lambda^\perp$ is {\it spherical\/} in $M^n$, that is, 
the leaves of $E_\lambda^\perp$ are umbilical submanifolds of $M^n$ with parallel mean curvature vector.

  It was pointed out in recent work due to Cecil and Jensen  
that there are two natural settings  for attempting to 
obtain classification results for proper Dupin hypersurfaces. 
One can either assume compactness and look for global results 
or work locally and  search for hypersurfaces that are locally irreducible. A Dupin hypersurface is {\it reducible\/} 
if it is Lie equivalent to a hypersurface obtained by one of Pinkall's constructions. 
Cecil and Jensen (\cite{cj}) showed that a locally irreducible
proper Dupin hypersurface with three distinct principal curvatures  must be Lie equivalent to an isoparametric hypersurface. 
On the other hand, Pinkall and Thorbergsson \cite{pt}, and independently Miyaoka and Ozawa \cite{mo}, produced compact 
proper Dupin hypersurfaces with $4$ distinct principal 
curvature that are neither locally Lie equivalent to 
isoparametric hypersurfaces nor locally reducible.
Thus, classifying locally or globally proper Dupin 
hypersurfaces with at least $4$ distinct principal 
curvatures remains wide open and seems to be a rather 
difficult problem.

The results in this article give strong support to our belief that a weaker notion of reducibility is more appropriate for 
the local study of Dupin hypersurfaces with arbitrary number of principal curvatures. We say that a proper Dupin hypersurface \mbox{$f\colon\,M^n\to\R^{n+1}$} is {\it weakly reducible\/} if it has a principal curvature $\lambda$ with integrable conullity $E_\lambda^\perp$, a property that is also invariant under Lie transformations. As observed before, every reducible Dupin hypersurface is also weakly reducible, but we will show that the converse does not hold. Another important subclass of weakly reducible Dupin hypersurfaces is that of  the
{\em holonomic\/} ones, that is, hypersurfaces that can be endowed with principal coordinates. In fact, holonomicity can be characterized by the fact that the submanifold  is weakly reducible with respect to every 
principal curvature. 
On the other hand, no isoparametric hypersurface of the sphere with at least three principal curvatures nor any of the examples in \cite{pt} and \cite{mo} are weakly reducible.

In this paper we give a complete local characterization of weakly reducible Dupin hypersurfaces. In fact, we solve a much more general problem with an interest of its own in the theory of Euclidean submanifolds of arbitrary codimension. Namely, we characterize the submanifolds that carry a Dupin principal normal with integrable conullity.

Recall that a smooth normal vector field $\eta$ of an isometric immersion \mbox{$f\colon\,M^n\to\R^{N}$} 
is called a  {\it principal normal\/} with multiplicity 
$s\geq 1$ if the tangent subspaces 
$$
\E_\eta=\ker (\a_f -\<\;\;,\;\;\>\eta)
$$ 
have constant dimension $s$, where 
$\a_f\colon TM\times TM\to T^{\perp}_fM$ stands for the second 
fundamental form of $f$ with values in the normal bundle. This is a natural generalization for submanifolds of higher codimension of the notion of principal curvature of a hypersurface. An important class of submanifolds that carry principal normals is that of submanifolds with flat normal bundle; cf.\ (\ref{nor2}). We say that a principal 
normal $\eta$ of multiplicity $s$ is {\em Dupin\/} if it is parallel in the normal connection of $f$ along the  (conformal) {\it  nullity \/} distribution $\E_\eta$ associated to $\eta$. This condition is automatic for multiplicity 
$s\geq 2$ (cf.\ \cite{re1} or \cite{dft}). If $\eta$ is  nonvanishing, it is  well-known  that  $\E_\eta$
is an  involutive distribution whose leaves are round \mbox{$s$-dimensional} spheres in $\R^N$; see \cite{re} or \cite{dft} for details. When $\eta$ vanishes identically, 
the distribution $\E_\eta=\E_0$ is known as the relative 
nullity distribution, in which case the leaves are open 
subsets of affine subspaces~of~$\R^N$.

      If one of the  first three constructions due to Pinkall 
is applied to an arbitrary submanifold  $L^{n-s}$ in 
$\R^{N-s}$ with any codimension $N-n$, then the process introduces a Dupin principal normal $\eta$ with multiplicity $s$, which has the additional property that the conullity $\E_\eta^\perp$ is a spherical distribution on $M^n$. 
It was proved in \cite{dft} 
(see also Theorem \ref{umb} below) that this last property characterizes these examples up to conformal transformations 
of the ambient space.

A simple way to construct submanifolds carrying a relative nullity distribution with integrable conullity is as follows.  
Let $g\colon\,L^{n-s}\to\Q_\e^{N}$, $\e=0,1,-1$, be an \ii with a parallel flat  normal subbundle  $\Nu$ of rank $s$.   
Then the $n$-dimensional {\it generalized cylinder\/} in $\Q_\e^{N}$ over $g$ determined by $\Nu$ is the  submanifold parametrized by means of the exponential map of  
$\Q_\e^{N}$ as
$$
\gamma\in\Nu\mapsto\mbox{exp}^\e_{g(\pi(\gamma))}(\gamma).
$$
Any such submanifold carries a relative nullity distribution of dimension $s$, whose leaves are the fibers of $\Nu$. 
Moreover, the conullity distribution is integrable, its leaves being given by the parallel sections of $\Nu$. 
Our first result concerning generalized cylinders is that  
these are the only submanifolds  having a relative nullity distribution with integrable conullity.

 The property of having a 
Dupin principal normal with integrable conullity is invariant under $\Lal$-transformations. By an {\it $\Lal$-transformation} of an Euclidean submanifold we mean a diffeomorphism that is a composition of conformal 
transformations of the  ambient space and parallel 
translations, the latter being translations of the submanifold 
by  parallel normal vector fields.  Of course, in the case of hypersurfaces $\Lal$-transformations are the usual  transformations of the Lie sphere geometry.
We call two submanifolds 
{\em $\Lal$-equivalent\/} if they differ by an \mbox{$\Lal$-transformation}. Therefore, a class of submanifolds carrying a Dupin principal normal with integrable conullity is obtained by applying $\Lal$-transformations to the family of (stereographic projections of) generalized cylinders. In the hypersurface case, this class properly contains those submanifolds obtained by Pinkall's constructions. However, they are far from exhausting the whole family of submanifolds carrying  a Dupin principal normal with integrable conullity,  as will be made clear below.

      The key observation in the characterization of submanifolds carrying a Dupin principal normal with integrable conullity is that the leaves of the conullity distribution of such a submanifold are always Ribaucour transforms one 
of each other. This is in the sense of the extended notion of Ribaucour transformation for submanifolds of arbitrary 
dimension and codimension developed in \cite{dt1} and 
\cite{dt2} from the classical notion for surfaces in three dimensional Euclidean space. This observation can be seen as a generalization of the classical fact (see \cite{bi}) that the orthogonal surfaces 
of a cyclic system are Ribaucour transforms one of each other. It has also been made recently by Corro \cite{co} in the particular case of holonomic Dupin hypersurfaces with a principal curvature of constant multiplicity one.  
       
       In order to turn the above observation into an explicit description of all such submanifolds, it was convenient to introduce the notion of \mbox{$\N$--Ribaucour} transform of a submanifold $h\colon\;L^{n-s}\to \R^{N}$ carrying a parallel flat normal subbundle $\N$ of rank $s$. This is an explicitly parametrized $n$-dimensional submanifold foliated 
by Ribaucour transforms of $h$, each one corresponding 
to a parallel section of $\N$. The orthogonal distribution 
to this foliation is precisely the nullity distribution of a Dupin principal normal.  One of the main results of this paper is that any submanifold that carries a Dupin principal normal with integrable conullity arises locally this way. 
       
  Each $\N$-Ribaucour transform of $h\colon\,L^{n-s}\to \R^{N}$ 
is essentially determined by a Codazzi tensor on $L^{n-s}$ that commutes with the \sff of $h$. 
We show that submanifolds that are $\Lal$-equivalent to 
(stereographic projections of) generalized cylinders are precisely those $\N$--Ribaucour transforms of $h$ that are determined by commuting Codazzi tensors 
on $L^{n-s}$ that can be expressed 
as linear combinations of the identity tensor and shape operators with respect to parallel normal vector fields. 

The results discussed in the preceding paragraphs are then 
applied to the class of \mbox{$k$--Dupin} submanifolds,  
that is, Euclidean submanifolds with flat normal \mbox{bundle} 
that have exactly $k$ principal normals all 
of which are Dupin.    
Our main result is that any \mbox{$k$--Dupin} submanifold that is {\it weakly reducible}, that is, carries a principal normal with integrable conullity, is  
the $\N$-Ribaucour transform of a $(k\!-\!1)$--Dupin
submanifold determined by a commuting \mbox{Codazzi}
tensor of Dupin type. 
A major application of this result is obtained by applying 
it to the important subclass of holonomic $k$--Dupin submanifolds. Namely, an explicit recursive procedure 
is derived to generate all such submanifolds in terms 
of~solutions of completely integrable systems of linear 
partial differential equations
of first order. 

We conclude the paper with  several additional results on
$k$--Dupin submanifolds. 
We show that the maximum possible value for the conformal codimension is $k-1$.
Moreover, the submanifold is necessarily holonomic if its conformal codimension is $k-1$ and it is necessarily weakly reducible if its conformal codimension is at least \mbox{$(2/3)k-1$}, the latter estimate being sharp.
Finally, we give a complete description of the 
weakly reducible \mbox{$4$--Dupin}  submanifolds.
We show that the submanifold is either holonomic or it is \mbox{$\Lal$-equivalent} to (the stereographic projection of) a generalized cylinder over a submanifold that is Lie 
equivalent to an isoparametric hypersurface. 
\medskip 

We are very grateful to T. Cecil and C. Olmos for several
helpful comments. \vspace{3ex}

\noindent{\bf\large  The Ribaucour transformation} 
\vspace{3ex}

In this section, we first recall the notions of Combescure 
and Ribaucour transformation of an Euclidean submanifold. 
Then, we discuss several basic 
facts about them that are used throughout the paper. 
We refer to \cite{dt1} and \cite{dt2} for further details and results on the subject.  
\vspace{1ex}
 
 A smooth map $\Fes\colon\,M^n\to\R^{n+p}$ is said to be a 
{\em Combescure transform\/} of a given \ii 
$f\colon\, M^n\to\R^{n+p}$ when there exists a symmetric endomorphism $\Phi$ of $TM$ such that 
$$
\Fes_*=f_*\circ \Phi.
$$ 
This condition  implies that $\Phi$ belongs  to the real 
vector space of (symmetric) Codazzi tensors on $M^n$ that 
are {\it commuting} 
in the sense that 
$$
\a_f(X,\Phi Y)=\a_f(\Phi X,Y) \fall X,Y\in TM.
$$
Conversely, any commuting Codazzi tensor $\Phi$ on a simply connected 
$M^n$ gives rise to a Combescure transform $\Fes$ of $f$. Moreover, 
$\Phi$ and $\Fes$ can be given as 
$$
\Phi=\hess\varphi - A^f_\beta \;\;\;\;\mbox{and}\;\;\;\; \Fes=f_*\nabla\va +\beta,
$$
where $\va\in C^{\infty}(M)$ and $\beta\in T_f^{\perp}M$ 
satisfy
\be\label{eq:gnorm}
\a_f(\nabla \va,X)+\nap_X \beta=0 \fall X\in TM,
\ee
$A^f_\beta$ denotes the shape operator of $f$ with 
respect to $\beta$ and $\nap$ stands for the induced 
connection  on the normal bundle. 
\vspace{1.5ex}

\noindent{\bf Definition.\/}  An immersion
$\tilde{f}\colon\, M^n\to\R^{n+p}$ is called a 
{\em Ribaucour transform\/} 
of a given immersion $f\colon\, M^n\to\R^{n+p}$ 
if $\tilde{f}\neq f$ 
everywhere, and there are a vector bundle isometry $\P\colon f^*T\R^{n+p}\to\tilde{f}^*T\R^{n+p}$ with 
$\P\,T_fM=T_{\tilde f}M$, and a nowhere vanishing smooth map 
$\delta\colon\, M^n\to\R^{n+p}$ such that:
\begin{itemize}
\item[(a)] $\P Z-Z=\<\delta,Z\>(f-\tilde{f})
\fall Z\in f^*T\R^{n+p}$;
\item[(b)] $D=f_*^{-1}\P^{-1}\tilde f_*\colon\,TM \to TM$ is 
self adjoint in the metric induced by $f$.
\end{itemize}\vspace{1ex}
 
 Condition $(a)$ says
that for any $Z\in T_{f(x)}\R^{n+p}$ the straight lines in $\R^{n+p}$ 
through $f(x)$ and $\tilde{f}(x)$ tangent to $Z$ and 
$\P Z$,\, respectively,\, are either parallel or intersect 
at a point 
equidistant to $f(x)$ and~$\tilde{f}(x)$. 
\vspace{1ex}

  The following statement contains the basic facts on the Ribaucour transformation that will be used throughout this 
paper without further reference. 
We denote by $\Sal(f)$ the set of pairs $(\va,\beta)$ 
satisfying (\ref{eq:gnorm}) such that 
$\va\Fes\neq 0$ 
everywhere. Then $\Sal_0=\Sal_0(f)$ stands for the corresponding real projective set, that is, $(\varphi, \beta)\sim(\varphi', \beta')$ if and only if 
$\varphi'= \lambda\varphi$ and 
$\beta' = \lambda\beta$ for some $0\neq\lambda\in\R$.

\begin{theorem}{\em (\cite{dt2})}{\bf .}\label{cod} Let $f\colon\,M^n\to\R^{n+p}$ be an \ii of a simply connected Riemannian manifold and let $\tilde f$ be a  Ribaucour 
transform of $f$. 
Then there exists a unique $[(\varphi,\beta)]\in \Sal_0$ such~that
\be\label{eq:rb}
\tilde{f} = f - 2\va\nu\Fes,
\ee
where $\Fes=f_*\nabla\va +\beta$ and $\nu=\|\Fes\|^{-2}$.  Moreover, 
we have that
$$
\P=I - 2\nu \Fes\Fes^*, \;\;\; \delta =-\varphi^{-1}\Fes 
\;\;\;\;\mbox{and}\;\;\;\; D=I-2\va\nu\Phi,
$$
where $\Fes^*(Z)=\<\Fes,Z\>$ for any $Z\in f^*T\R^{n+p}$.
Conversely, given $[(\varphi,\beta)]\in \Sal_0$, let 
$\P$, $\delta$ and $D$ be defined by the preceding expressions 
on an open subset 
$U\subset M^n$ where $D$ is invertible. Then 
$\tilde{f}|_{U}$ given by {\em (\ref{eq:rb})} 
is a Ribaucour transform of $f|_{U}$ for $\P$, $\delta$ 
and $D$. Moreover:\vspace{1ex}
 
\noindent $(i)$ The second fundamental forms of $f$ and 
$\tilde f$ are related by 
$$
A^{\tilde f}_{\P \xi}=D^{-1}(A^f_\xi + 2\nu\<\beta,\xi\>\Phi);
$$

\noindent $(ii)$ The restriction 
$\P|_{T_f^\perp M}\colon\,T_f^\perp M\to T_{\tilde f}^\perp M$ 
is parallel.
\end{theorem}
 
  We denote by $\Ral_w(f)$ the Ribaucour transform $\tilde{f}$ 
of $f$ determined by $w\in \Sal_0$. Since $w$ determines  $\P,\delta$ and  $D$ completely and $\Fes$ and $\Phi$ up to  constants, when convenient we will use it as a subscript for these maps.\vspace{1ex}

   We see next that inversions and parallel translations are special cases of Ribaucour transformations. In the following 
and elsewhere writing a vector subspace as a subscript of a vector indicates taking the orthogonal projection of the vector onto that subspace.

\begin{examples}\po\label{ribs}{\em
\noindent $(i)$ Given a point $P_0\in \R^{n+p}$ and $r>0$, set $w=[(\va_1,\beta_1)]$ where
$2\va_1=\|f-P_0\|^2-r^2$ and $\beta_1=(f-P_0)_{T_f^{\perp}M}$.  
Then $\Fes=f-P_0$, $\Phi=I$, and
$$
\tilde f=\Ral_{w}(f)=P_0+r^2\|f-P_0\|^{-2}(f-P_0)
$$ 
is obtained from $f$ by an inversion with respect to 
the sphere of radius $r$ centered at $P_0$.  Moreover, 
$\P=\id -2\|f-P_0\|^{-2}(f-P_0)^*(f-P_0)$ and
\be\label{sff}
r^2A^{\tilde f}_{\P\mu}=\|f-P_0\|^2A^f_\mu+2\<f-P_0,\mu\>\id
\fall\mu\in T_f^\perp M.
\ee

\noindent $(ii)$ Given a parallel normal vector field $\xi$, 
set \mbox{$w=[(\va_2,\beta_2)]$} 
where $2\va_2=\|\xi\|^2$ and $\beta_2=-\xi$. 
Then $\Fes=-\xi$, $\Phi=A^f_{\xi}$, and 
$$
\tilde f=\Ral_{w}(f)=f+\xi
$$ 
is the parallel translation $L_\xi$ of $f$. Moreover, $\P=\id-2\|\xi\|^{-2}\xi^*\xi$ and
\be\label{sff2} 
A^{\tilde f}_\mu=(\id-A^f_\xi)^{-1}A^f_\mu 
\fall \mu\in T_f^\perp M.
\ee
}\end{examples}

   The Ribaucour transformation has the following invariance property  under $\Lal$- transformations (see Proposition $33$ in \cite{dt2}).  By an {\it $\Lal$-transformation} of an Euclidean submanifold we mean a diffeomorphism that is a composition of conformal transformations of the  ambient space and parallel translations, the latter being translations of the submanifold by  parallel normal vector fields. To each 
$\Lal$-transformation $T$ of $f$ we associate an~$\Lal$-transformation $\tilde T$ of 
a given Ribaucour transform $\Ral_w(f)$ as follows: 
\begin{itemize}
\item[(i)] $\tilde T=T$ when $T$ is a conformal map 
of $\R^{n+p}$; 
\item[(ii)] $\tilde T=L_{\P\xi}$ when $T=L_\xi$.
\end{itemize}
Then, there is a correspondence 
$T\longmapsto w^T\in \Sal_0(T(f))$
such that
\be\label{ruso}
\tilde T(\Ral_w(f)) =\Ral_{w^T}(T(f)).
\ee
For later use we describe  that correspondence 
explicitly for each of the following types of transformations. We omit the computations, which are straightforward with the exception of  case $(iv)$. For the latter we refer to Proposition~$31$ in \cite{dt2}.
\vspace{1ex}

\noindent $(i)$ {\it Euclidean translation:} $\;T_u(f)=f+u,     
\;\mbox{where}\;\; u\in\R^{n+p}$. Then,
$$
\Fes^{T_u} = \Fes\;\;\mbox{and}\;\; \varphi^{T_u}=\varphi.
$$

\noindent $(ii)$ {\it Orthogonal transformation:} 
$\;O(f)=O\circ f,\;\mbox{where}\;\; O \in O(n+p)$. Then,
$$
\Fes^{O}=O(\Fes)\;\;\mbox{and}\;\; \varphi^{O}=\varphi.
$$

\noindent $(iii)$ {\it Homothety:} $\;H_k(f)=kf,\;\mbox{where}\;\; k\in \R$. Then,
$$
\Fes^{H_k}=\Fes\;\;\mbox{and}\;\;\varphi^{H_k}=k\varphi.
$$

\noindent $(iv)$ {\it Inversion:} $\;i(f)=f/\|f\|^2$. Then,
$$
\Fes^{i}=\P^i(\Fes-2\varphi \|f\|^{-2}f)
=\Fes - 2(\<\Fes,f\>-\varphi)\|f\|^{-2}f\;\;\mbox{and}\;\; 
\varphi^{i}=\varphi\|f\|^{-2}.
$$

\noindent $(v)$ {\it Parallel translation:} 
$\;L_\xi(f)= f+\xi$, where $\xi\in T_f^\perp M$ is parallel. 
Then,
$$
\Fes^{L_\xi}=\Fes \;\;\mbox{and} \;\;\varphi^{L_\xi}=\varphi+\<\Fes,\xi\>.
$$

The {\it conformal codimension\/} of a submanifold $g\colon\,M^n\to 
\R^{n+p}$ is the number
$c(g)$  such that $n+c(g)$  is the least dimension of a sphere or an 
affine subspace in $\R^{n+p}$ that contains
the submanifold. If $c(g)=p$ then $g$ is said to be {\it conformally 
substantial\/}.

\begin{proposition}\po\label{inv} The conformal codimension is invariant 
under $\Lal$-transformations.
\end{proposition}

\proof The invariance under conformal transformations is clear, thus it 
suffices to check that a submanifold $g\colon\,M^n\to \R^{n+p}$ and a 
parallel translate $L_\xi(g)$ have the same conformal codimension. 
Since the immersions have the same normal spaces at every point, it 
suffices to argue that  $L_\xi(g)(M)$ is contained in a sphere whenever $g(M)$ is 
contained in a sphere centered, say, at the origin. We can write $\xi=ag+b\eta$, where 
$a,b\in \R$ and $\eta$ is a parallel normal vector field tangent to the 
sphere. Hence $g+\xi=(1+a)g+b\eta$ has also constant norm.\qed\vspace{1ex}

On a submanifold $f\colon\,M^n\to\R^{n+p}$ consider a commuting Codazzi 
tensor
\be\label{tensor}
\Phi=a\id-A^f_\delta,\;\;\; a\in\R,
\ee
where $\delta\in T_f^\perp M$ is a parallel vector field in the normal 
connection. Then, a Combescure transform $\Fes$ of $f$ determined by 
$\Phi$, that is,  $\Fes_*=f_*\circ \Phi$, can  be written as
\be\label{Fes}
\Fes = af + v +\delta,\;\;\; v\in\R^{n+p}.
\ee

\begin{proposition}\po\label{exp}   The decomposition
{\em(\ref{Fes})} is unique if $f$ is conformally substantial.
\end{proposition}

\proof  If $\Fes=\tilde a h + \tilde v + \tilde \delta$, set
$a'=\tilde a-a$, $v'=\tilde v-v$ and $\delta'=\tilde\delta-\delta$. We 
obtain differentiating $a'h+v'+\delta'=0$
that $A^h_{\delta'}=a'\id$. Since $\|\delta'\|$ is constant
and $h$ is conformally substantial, we conclude that
$\delta'= 0=a'=v'$.\qed\vspace{1ex}

If $w=[(\va,\beta)]\in\Sal_0(f)$ is such that 
$\Fes=\Fes_w=f_*\nabla\va+\beta$ is as in (\ref{Fes}), then
\be\label{fis}
2\varphi= a\|f\|^2 + 2\<f,v\> + c,\;\;\; c\in\R.
\ee

\noindent{\bf Definitions.\/} We say that
$w, \Phi$ or $\Fes$ are {\it $\Lal$-trivial\/} when they are given by  (\ref{tensor}), (\ref{Fes}) and  (\ref{fis}). 
They are {\it conformally trivial\/} if they can be given by those expressions  with $\delta=0$.
\vspace{1ex}

   If  $w\in\Sal_0(f)$ is $\Lal$-trivial, then the 
corresponding  Ribaucour transform is   
\be\label{gen}
\tilde{f} = f - \left(a\|f\|^2 + 2\<f,v\> + c\right)
\frac{af + v +\delta}{\|af + v +\delta \|^2}.
\ee
Notice that special examples in (\ref{gen}) are inversions and parallel translations.

\begin{proposition}\po\label{trivialt}
If $w\in\Sal_0(f)$ is $\Lal$-trivial (respectively, conformally trivial) then the same holds for $w^T\in\Sal_0(T(f))$ for any $\Lal$-transformation $T$ of $f$. More precisely, if $\Fes$ and $\va$ are given by {\em(\ref{Fes})} and {\em(\ref{fis})},  then $\Fes^T$ and $\va^T$ are as follows:

\begin{itemize} 
\item[{\em (i)}]
$\left\{\begin{array}{l}
\Fes^{T_u}=aT_uf + v-au + \delta,\vspace*{1.5ex}\\
2\varphi^{T_u}=a\|T_uf\|^2 + 
2\<T_uf,v-au\>+c-2\<u,v\>+a\|u\|^2,
\end{array}\right.$
\medskip
\item[{\em (ii)}]
$\Fes^{O}=aOf + Ov + O\delta, \;\;\;
2\varphi^{O}=
a\|Of\|^2 + 2\<Of,Ov\> + c$,
\medskip
\item[{\em (iii)}]
$\Fes^{H_k}=(a/k)kf + v + \delta, \;\;\;2\varphi^{H_k}=
(a/k)\|kf\|^2 + 2\<kf,v\>+ ck$,
\medskip
\item[{\em (iv)}]
$\Fes^{i}=ci(f) + v + \delta-2\<\delta,f\>i(f),
\;\;\;2\varphi^{i}=
c\|i(f)\|^2 + 2\<i(f),v\> + a$,
\medskip
\item[{\em (v)}]
$\left\{ \begin{array}{l}
\Fes^{L_\xi}=aL_\xi f + v + \delta-a\xi,\vspace*{1.5ex}\\
2\varphi^{L_\xi}=a\|L_\xi f\|^2 + 
2\<L_\xi f,v\> + c +2\<\delta,\xi\>-a\|\xi\|^2.
\end{array} \right.$
\end{itemize}
\end{proposition}

\proof It is straightforward using the expressions after (\ref{ruso}).\qed\vspace{1.5ex}

\vspace{2ex}

\noindent{\bf\large  Dupin principal normals} 
\vspace{3ex}

 Our goal in this section is to describe a procedure 
to generate all Euclidean submanifolds carrying a Dupin 
principal normal with integrable conullity. \vspace{1ex} 
  
 We start by observing that the property of carrying a Dupin principal normal with integrable conullity is  invariant under 
${\cal L}$-transformations. In fact, let  $f\colon\, M^n\to\R^{n+p}$ be an \ii and let $\tilde f$ be an ${\cal L}$-transform  of $f$. Then, to each principal normal $\eta$ 
of $f$ there corresponds  a principal normal $\tilde\eta$ of 
$\tilde f$ such that $\E_{\tilde\eta}=\E_\eta$, and thus $\tilde\eta$  has integrable conullity if and only if $\eta$ does. Moreover, $\tilde\eta$ is Dupin if and only if the same holds for $\eta$. Namely,
$$
\tilde\eta=r^{-2}\P(\|f-P_0\|^2\eta+2(f-P_0)_{T_f^\perp M}),
$$
for an inversion as in Examples \ref{ribs} - $(i)$, whereas
$$\tilde\eta=(1-\<\xi,\eta\>)^{-1}\eta
=(1-\<\xi,\eta\>)^{-1}\P(\eta-2\|\xi\|^{-2}\<\eta,\xi\>)$$
for a parallel translation $L_\xi$ as 
in Examples \ref{ribs} - $(ii)$, as one can easily check  
using (\ref{sff}) and (\ref{sff2}).
That $\tilde\eta$ is Dupin if and only if the same is true for $\eta$  follows  from  $(ii)$ in Theorem~\ref{cod}\vspace{1ex}.

   We now introduce the main tool of the paper, which is an extension of the notion of Ribaucour transformation of an Euclidean submanifold.
 Fix a simply connected submanifold $h\colon L^{n-s}\to\R^{n+p}$ with 
a normal subbundle $\N$ of rank $s$  that is parallel and flat with respect 
to the normal connection. Flatness of $\N$ means that the normal curvature tensor satisfies $R_h^\perp|_{\N}=0$.
We denote by $\Sal_\N=\Sal_\N(h)$ the set of equivalence classes of pairs $(\va,\beta)\in \Sal(h)$  under 
the equivalence relation that identifies two pairs $(\varphi, \beta)$ and $(\varphi', \beta')$ whenever 
$$
\varphi'= \lambda\varphi\;\;\; \mbox{and}\;\;\; 
\beta' -\lambda\beta \in\Np
$$
for some $0\neq\lambda\in\R$.  Here and elsewhere $\Np$ stands for the $s$-dimensional real vector space of parallel sections of $\N$.\vspace{1ex}

\noindent{\bf Definition.\/} The 
{\it $\N$-Ribaucour transform\/} 
$\Ral^\N_w(h)$ of $h\colon L^{n-s}\to\R^{n+p}$ determined by $w=[(\va,\beta)]\in\Sal_\N$ is the $n$--dimensional immersed submanifold  parametrized, at regular points, 
by the map $f\colon L^{n-s}\times \Np\to\R^{n+p}$ given by 
\be\label{rib2}
f(u,t)=h_t(u),
\ee
where $h_t=\Ral_{[(\va,\beta+t)]}(h)
=h-2\va\nu(h_*\nabla\va+\beta+t)$ and $\nu=\|h_*\nabla\va+\beta+t\|^{-2}$.
\vspace{1ex}
 
 Observe that $h$ itself is the leaf ``at infinity" of the foliation parametrized by $t$ in the sense that $h=\lim_{\|t\|\to\infty}h_t$.  
After choosing a parallel orthonormal frame 
$\xi_1,\ldots,\xi_s$ of $\N$, the map $f$ can be rewritten 
as $f\colon L^{n-s}\times \R^s\to \R^{n+p}$ given~by
$$
f=h-2\va\nu(h_*\nabla\va+\beta+\sum_{j=1}^sy_j\xi_j(u)).
$$

Given $(\va,\beta)\in\Sal(h)$ and $t\in\Np$, we denote 
$\Fes_t=h_*\nabla\va+\beta+t$, $\nu_t=\|\Fes_t\|^{-2}$, $\Phi_t=\hess\va - A^f_{\beta+t}$,
$\P_t=\P_{[(\va,\beta+t)]}$ and $D_t=D_{[(\va,\beta+t)]}$.

\begin{proposition}\po\label{multi} Let  
$\Ral^\N_w(h)\colon\,M^n=L^{n-s}\times \Np\to\R^{n+p}$ be an $\N$-Ribaucour transform of 
$h\colon L^{n-s}\to\R^{n+p}$. Then the following facts hold:
\vspace{1ex}

\noindent $(i)$ $f=\Ral^\N_w(h)=\Ral^\N_{\bar w}(h)$ if and 
only if $w=\bar w$.\vspace{1ex}

\noindent $(ii)$ The normal space  of $f$ at $(u,t)\in M^n$
is $T^\perp_{(u,t)}M=\P_t(\N^\perp(u))$.  
\vspace{1ex}
 
\noindent $(iii)$ The normal connection of $f$ is given by
$$
\nap_X\hat\delta=\P_t\nap_X\delta\;\;\;\mbox{and}\;\;\;
\nap_S\hat\delta=0 
$$
for all $\delta\in \N^\perp, X\in TL$ and $S\in\Np$, where $\hat\delta(u,t) 
= \delta_t(u) = \P_t\delta(u)$. In particular,
$\delta$ is a parallel normal  vector field of $h$ if and 
only if $\hat\delta$ is a parallel normal vector field of $f$, and hence, $f$ has flat normal bundle if and only if $h$ 
does.\vspace{1ex}

\noindent $(iv)$ The second fundamental form of $f$ 
at $(u,t)\in M^n$ 
for all $S\in\Np$ and $Z\in T_{(u,t)}M$ is given by
\be\label{una}
\alpha_f(S,Z)=\<S,Z\>\P_t\bar{\beta}(u),
\ee
where $\bar\beta=-\va^{-1}\beta_{\N^\perp}$, and  by
\be\label{dos}
\alpha_f(X,Y)=\P_t((\alpha_h(D_tX,Y)
+2\nu_t(u)\<D_tX,\Phi_tY\>\beta(u))_{\N^\perp})
\ee
for all $X,Y\in T_uL$.

\noindent $(v)$ For each point $u_0\in L^{n-s}$ the map $f(u_0,\,\cdot\,)\colon\,\Np\to\R^{n+p}$  
is a conformal parametrization of a sphere or an
affine subspace, the latter occurring if and only 
if $\Fes_{\va,\beta}(u_0)\in\N$.
\end{proposition}

\proof The proof of $(i)$ is straightforward. An easy 
computation at $(\,\cdot\,, t)$ yields
\be\label{tang}
f_*X=\P_tD_tX \;\;\;\mbox{and}\;\;\;
f_*S=-2\va\nu_t\P_tS 
\ee
for all $X\in TL$ and $S\in\Np$, and $(ii)$ follows. The assertions in $(iii)$ are consequences of $(ii)$ in
Theorem \ref{cod} and the fact that  
$\nab_S\hat\delta=-2\nu_t\<\beta,\delta\>\P_tS\in TM$,
where $\nab$ denotes the Euclidean connection.
The proof of $(iv)$ is similar to that of Corollary $21$ of \cite{dt2}. Finally, $(v)$ follows  from  $(\ref{una})$ and (\ref{tang}).\vspace{1.5ex}\qed

     We show next that the $\N$-Ribaucour transform  
$\Ral^\N_w(h)$ of an \ii $h\colon L^{n-s}\to\R^{n+p}$ 
with a parallel flat normal subbundle $\N$ of rank $s$ 
always carries a Dupin principal normal of multiplicity $s$ 
if a certain regularity assumption is satisfied. For the 
precise statement  we need to introduce some terminology. 
For an \ii $h\colon\;L^{n-s}\to \R^N$, a  vector subbundle $\Nu\subset T^\perp_hL$ and a vector field $\delta\in \Nu$, 
we define at $x\in L^{n-s}$ the subspace
$$
\E^\Nu_\delta(x)=\{Z\in T_xL :\,\a_h(Z,X)_\Nu=\<Z,X\>\delta
\fall X\in T_xL\}.
$$
Given $w=[(\va,\beta)]\in\Sal_\N(h)$, we denote
$\bar\beta=-\va^{-1}\beta_{\N^\perp}$
and $\E(w)=\E^{\N^\perp}_{\bar\beta}$.
\vspace{1ex}

The only purpose of the regularity condition 
$\E(w)=0$ in the following result is to assure that the  
Dupin principal normal generated by the \mbox{$\N$-Ribaucour} transformation has the lowest possible multiplicity everywhere.
  
\begin{proposition}\po\label{1} 
Let $h\colon L^{n-s}\to \R^{n+p}$ 
be a simply connected submanifold carrying a 
parallel flat normal subbundle 
$\N$ of rank~$s$, and let $$
f=\Ral^\N_w(h)\colon\,M^n=L^{n-s}\times \Np\to\R^{n+p}
$$ 
be an $\N$-Ribaucour transform of $h$ determined by 
$w\in\Sal_\N$ with $\E(w)=0$. Then the vector field 
$\eta(u,t)\!=\!\P_t\bar\beta(u)$ is a Dupin 
principal~normal of $f$ with integrable conullity, the leaves being parametrized by \mbox{$(\,\cdot\,,t_0)$ with $t_0\in\Np$}.
\end{proposition}

\proof It follows from (\ref{una}) that $\Np\subset\E_\eta$.
To see that equality holds, take $X\in TL\cap\E_\eta$ and 
$Y\in TL$. Then
$$
\a_f(X,Y)=\<D_tX,D_tY\>\P_t\bar\beta.
$$
Using (\ref{dos}),  we obtain that
$$
\a_h(D_tX,Y)_{\N^\perp}=\<D_tX,Y\>\bar\beta,
$$
and hence $X=0$ by the assumption that $\E(w)=0$.  That
$\eta$ is Dupin follows from Proposition \ref{multi} - $(iii)$. Finally, by (\ref{tang}) we have that  $\<f_*S,f_*X\>=0$, 
and we conclude that the conullity is integrable, the leaves being parametrized by~\mbox{$(\,\cdot\,,t_0)$ with $t_0\in\Np$}.\vspace{1.5ex}\qed

One of our main results is that, conversely, any submanifold carrying a Dupin principal  normal with integrable conullity arises locally this way.

\begin{theorem}\po\label{thm:main} Let 
$f\colon M^n\to \R^{n+p}$ be an \ii carrying a Dupin principal  normal $\eta$ of multiplicity $s$ with integrable conullity. 
Then $\N=\E_\eta|_{L}$ is a parallel flat normal subbundle of 
$h=f|_{L}$ for any given leaf $L^{n-s}$ of conullity, and $f$ 
is an \mbox{$\N$-Ribaucour} transform of $h$ determined by a unique $w\in\Sal_\N(h)$ with $\E(w)=0$.
\end{theorem}
 
      Theorem \ref{thm:main} will be derived from Proposition \ref{marc},
where the conullity is only assumed to admit one {\it maximal integral submanifold\/}, that is, a submanifold $L$ satisfying that $T_xL=\E_\eta^\perp(x)$ for any $x\in L$, and that 
each leaf of $\Delta_f$ intersects $L$ exactly once.
First we need some further terminology and notations. 
We refer to a vector bundle $(E,\pi, M)$ with total 
space $E$ and 
projection \mbox{$\pi\colon\,E\to M$} simply by $E$, and 
denote by $\Gamma(E)$ 
the space of its smooth sections. The kernel of $\pi_{*t}$ at $t\in E(x)$ 
is the vertical subspace of $T_tE$. Clearly, $T_tE$
can be identified with $E(x)$ itself, thus we can write $E(x)\subset T_tE$ 
without risk of confusion.  We call the {\it vertical subbundle\/} 
of $TE$ the subbundle $E\subset TE$ formed by its 
vertical subspaces. 

\begin{proposition}\po\label{marc}
Let $h\colon\, L^{n-s}\to\R^{n+p}$ be a submanifold with a parallel normal subbundle $\N$ of rank $s$ and let 
$\mu\colon\,L^{n-s}\to\R^{n+p}$ be a smooth map such that
\mbox{$\eta_h=\mu_{\N^\perp}$}  satisfies   
$\E_{\eta_h}^{\N^\perp}=0$ and
\be\label{nocombe}
(\mu_* X)_{\N^\perp}=\<h_*X,\mu\>\eta_h \fall X\in TL.
\ee
Define $f\colon\,\N\to\R^{n+p}$ by
\be\label{rib}
f\circ t = h + 2\|\mu+t\|^{-2}(\mu+t) \fall t\in\Gamma(\N),
\ee
and $\eta\colon\,\N\to\R^{n+p}$ by
\be\label{ribb}
\eta\circ t= \eta_h-2\|\eta_h\|^2\|\mu+t\|^{-2}(\mu+t).
\ee
Then $f$ parametrizes, at regular points, a submanifold $M^n$ with $\eta$ 
as a Dupin principal  normal such that the nullity $\E_\eta$ 
is the vertical subbundle $\N\subset T\N$, and $L^{n-s}$ is a
maximal integral submanifold of the  conullity.

Conversely, let $f\colon\,M^n\to\R^{n+p}$ be a submanifold carrying a 
Dupin principal  normal $\eta$ of multiplicity $s\geq 1$ such  that the conullity has a maximal integral submanifold $L^{n-s}$. Then 
$\N=\E_\eta|_{L}$ is normal  and parallel along $h=f|_{L}$, and there 
exists a smooth map $\mu\colon\, L^{n-s}\to\R^{n+p}$ orthogonal to  $\N$ 
such that $f$  and $\eta$   can be locally parametrized  by 
{\em (\ref{rib})} and {\em (\ref{ribb})} on an open 
neighborhood of $L^{n-s}$.
\end{proposition}

     Before proving Proposition \ref{marc} we make the 
following useful observation.

\begin{lemma}\po\label{back} Let $h\colon\, L^{n-s}\to\R^{n+p}$ be a submanifold with a parallel normal subbundle $\N$ of rank $s\geq 0$ and let $\eta \in \N^\perp$ be a nowhere vanishing vector field. Then the subspace 
$\{T\in \E^{\N^\perp}_\eta(x):\nabla_T^\perp \eta=0\}$ 
coincides with 
$$
\{T\in T_xM :(h+\|\eta\|^{-2}\eta)_*T=0;\; \tilde{\nabla}_T \xi\in\{\eta\}^\perp\cap \N^\perp\;\;\mbox{{\em for all}}\;\;\xi\in \{\eta\}^\perp\cap \N^\perp\}.
$$
\end{lemma}

\proof We have that $T\in \E^{\N^\perp}_\eta(x)$ and $\nabla_T^\perp \eta=0$ if and only if the right hand sides of the following equations vanish for any 
$\xi\in \{\eta\}^\perp\cap \N^\perp$:
\begin{eqnarray*} (h+\|\eta\|^{-2}\eta)_*T \!\!&=&\!\! h_*(T-\|\eta\|^{-2}A_\eta T) +\nabla_T^\perp \|\eta\|^{-2}\eta;\\
\<\tilde{\nabla}_T\xi,\eta\> \!\!&=&\!\! -\<\nabla_T^\perp\eta,\xi\>;\\
\<\tilde{\nabla}_T\xi,h_*X\> \!\!&=&\!\! -\<\alpha_h(T,X),\xi\>.
\;\;\;\;\;\;\vspace{1ex}\qed
\end{eqnarray*}

\noindent{\em Proof of Proposition \ref{marc}:\/} 
We start the proof of the direct statement by 
determining the normal bundle of $f$. Differentiating 
(\ref{rib}) along $X\in TL$, we obtain that
\be\label{dife}
(f\circ t)_*X=h_*X+2X(\|\mu+t\|^{-2})(\mu+t) 
+ 2\|\mu+t\|^{-2}(\mu_*X+ \nab_Xt).
\ee
It follows from (\ref{nocombe}), (\ref{dife}) and the 
parallelism of $\N$
that $\N^\perp\cap\{\eta_h\}^\perp$ is normal to $f$. 
Moreover, (\ref{nocombe}) and (\ref{dife}) also imply that
$$
\<(f\circ t)_*X,\eta\circ t\>=2\|\mu+t\|^{-2}
\left(\<\mu_*X,\eta_h\>-\|\eta_h\|^2\<h_*X,\mu\>\right)=0.
$$
Using that $\|\eta\circ t\|=\|\eta_h\|$, we obtain the 
orthogonal splitting 
\be\label{net}
T_f^\perp M(t(x))
=\left(\N^\perp\cap\{\eta_h\}^\perp\right)(x)
\oplus\spa\{\eta(t(x)\},
\ee  
which shows that $\eta^\perp$ is constant along 
$\N\subset T\N$. For a point $x\in L^{n-s}$ where 
$\eta_h(x)=0$, we conclude from (\ref{net}) that
\be\label{nor}
T_f^\perp M(t(x))=\N^\perp(x),
\ee 
and thus $\N(x)\subset\E_0(t(x))$.  On the open subset where $\eta_h\neq 0$, we have that
\be\label{const}
f\circ t+\|\eta\circ t\|^{-2}\eta\circ t=
h+\|\eta_h\|^{-2}\eta_h,
\ee 
and hence the left hand side is constant along the leaves of $\N\subset T\N$.  We conclude from Lemma \ref{back} that $\N\subset\E_\eta$ and that $\eta$ is parallel along~$\N$.

 It remains to show that $\N=\E_\eta$ under our 
regularity assumption. This holds by (\ref{nor}) at the points
of $L^{n-s}$ where $\eta_h=0$. Hence, we may assume that $\eta_h(x)\neq 0$, 
and then the same holds for $\eta$ on $\N(x)$.
Any transversal tangent vector to the leaves of 
$\N\subset T\N$ can be written as~$t_*X$  for some 
$X\in T_xL$ and $t\in\Gamma(\N)$. Assume that $t_*X\in\E_\eta(t(x))$. Then, we have from Lemma \ref{back} 
that 
$$
\left((f+\|\eta\|^{-2}\eta)_*(t_*X)\right)_{T_fM}=0=
(\nab_{t_*X}\xi)_{T_fM}
$$ 
for any normal vector field $\xi\in\{\eta\}^\perp\subset T_f^\perp M$.
By (\ref{net}), any normal vector field \mbox{$\widehat{\xi}\in\N^\perp\cap\{\eta_h\}^\perp$} to $h$ 
gives rise to a normal vector field $\xi\in \{\eta\}^\perp$ to 
$f$ along $t$ by setting $\xi\circ t=\widehat{\xi}$. 
Therefore, $(\nab_{X}{\widehat{\xi}})_{T_hL}=0$
for any $\widehat{\xi}\in \N^\perp\cap\{\eta_h\}^\perp$. 
Moreover, by (\ref{const}) we obtain that
$\left((h+\|\eta_h\|^{-2}\eta_h)_*(X)\right)_{T_hL}=0$.
 From our assumption on $\eta_h$ and Lemma~\ref{back}, we conclude that $X=0$, and thus 
$t_*X=0$ as we wanted.
 
We first prove the converse under the assumption that $\eta$ 
is nowhere vanishing. Let $\sigma(x)$ be the leaf of $\E_\eta$ through 
$x\in L^{n-s}$. Then $f(\sigma(x))$ is an open 
subset of a round sphere $\Sp^s(x)$ through $f(x)$. Now let
\mbox{$\mu\colon\,L^{n-s}\to\R^{n+p}$} be the vector field defined so that $h(x)+\|\mu(x)\|^{-2}\mu(x)$ is the center of 
$\Sp^s(x)$. 
Then $\mu(x)$ is orthogonal to the tangent space $\N(x)=\E_\eta(x)$ 
of $\Sp^s(x)$ at $h(x)$, and the inversion with respect
to the sphere of radius $\sqrt{2}$ centered at the origin followed by translation by $h(x)$ maps the affine hyperplane 
$\mu(x)+\N(x)$ onto $\Sp^s(x)$ minus the point $h(x)$. 
Since $L$ is a maximal integral submanifold of the conullity 
then $f$ is parametrized by (\ref{rib}) on 
an open neighborhood of the zero section of 
$\N$  along $h$. Moreover, since  $\eta$ is a Dupin principal normal, it follows from Lemma \ref{back} that $f+\|\eta\|^{-2}\eta$ is constant along $\E_\eta$. Therefore (\ref{const}) 
holds and  $\eta$ is given by (\ref{ribb}). Finally, the 
proof of the direct
statement shows that $\eta$ being normal to $f$ and the 
normal subspace
$\{\eta\}^\perp$ being constant along $\E_\eta$ imply that 
$\mu$ satisfies (\ref{nocombe}) and that $\N$ is parallel.

  For the general case, we may compose $f$ with a translation, 
if necessary, and an inversion $i$  
so that the corresponding principal normal $\eta^i$ of $i(f)$ 
is nowhere vanishing. Thus, we have a submanifold $h$, a  subbundle $\N$ and a map $\mu$ as before, and we can describe $i(f)$ by (\ref{rib}) as 
$$
i(f)=h+2\|\mu +t\|^{-2}(\mu + t),
$$ 
with $\mu$ satisfying (\ref{nocombe}). 
Applying the inversion $i$ to $i(f)$, it is easy to see that
$$
f=i(h)+2\|\bar\mu+\bar t\|^{-2}(\bar \mu + \bar t),
$$
where $\bar\mu=\P_i(2h+\|h\|^2\mu)$, $\bar t=\P_i\,\|h\|^2t$, 
and  $\P_i=I-2\<h,\cdot\>i(h)$ is the vector bundle 
isometry associated to $i$ seen as a Ribaucour transformation 
of  $i(h)$. 
In particular, we have that $\bar\N^\perp=\P_i\,\N^\perp$.
It is also easy to check that 
$(\bar\mu_*X)_{\bar \N^\perp}=\<i(h)_*X,\bar \mu\>\bar\mu_{\bar\N^\perp}$,   
and this completes the proof.\qed\vspace{1.5ex}

\noindent {\em Proof of Theorem} \ref{thm:main}:  By Proposition \ref{marc}, 
we have that $\N=\E_\eta|_{L}$ is parallel along $h=f|_{L}$ for 
any given leaf $L^{n-s}$ of 
the conullity, and that $f$ and $\eta$ can be parametrized by (\ref{rib}) 
and (\ref{ribb}), respectively, for some smooth map 
\mbox{$\mu\colon\, L^{n-s}\to\R^{n+p}$} everywhere orthogonal 
to $\N$ satisfying (\ref{nocombe}).
Since the conullity is integrable, for any $T\in\N(x_0)$ 
there is a section $t\in \Gamma(\N)$ with
$t(x_0)=T$ everywhere orthogonal to the vertical subbundle of
$T\N$ with respect to the metric induced by $f$. We have that
$$
f_*S=2\|\mu+t\|^{-2}(S-2\|\mu+t\|^{-2}\<t,S\>(\mu+t))
$$
for any section $S$ of $\N\subset T\N$ along $t$. A straightforward 
computation using (\ref{dife}) shows that $\<(f\circ t)_*X,f_*S\>=0$ 
if and only if
$$
\<\mu_*X+\nab_X t,S\>-\<t,S\>\<h_*X,\mu\>=0.
$$
Since $\N$ is parallel this is equivalent to
\be\label{eq:t}
\nap_X t=-(\mu_*X)_\N+\<h_*X,\mu\>t.
\ee
Given $t_0\in \N(x_0)$,  take  $t_1, t_2\in \Gamma(\N)$ 
orthogonal to the vertical subbundle of
$T\N$ such that   $(t_1-t_2)(x_0)=t_0$. It follows from 
(\ref{eq:t}) that $t=t_1-t_2$  satisfies $t(x_0)=t_0$ and
\be\label{eq:nap}
\nap_X t = \o(X)t,\;\;\mbox{where}\;\; \o(X)=\<h_*X,\mu\>.
\ee
Set $\tau=\log (\|t\|/\|t_0\|)$.
Then $\tau(x_0)=0$ and $d\tau=\o$. It follows from 
(\ref{eq:nap}) and
closeness of $\o$ that $R^\perp(X,Y)t=0$ for all $X,Y\in TL$. 
In particular, we have that $\N$ is flat. Moreover, (\ref{nocombe}) and (\ref{eq:t}) yield
\be\label{eq:mu2}
(\mu_*X)_{T_h^\perp L}=X(\tau)(t+\mu_{\N^\perp})-
(\nab_Xt)_\N.
\ee
Set $t=e^\tau t_T+\bar{t}$, where $t_T$ is the parallel 
section of $\N$ 
with $t_T(x_0)=T=t_0$. Then (\ref{eq:mu2}) can be written as
$$
((\mu+\bar{t})_*X)_{T_h^\perp L}
=X(\tau)(\mu+\bar{t})_{T_h^\perp L},
$$
or equivalently,
$$
(e^{-\tau}(\mu+\bar{t}))_*X\in T_hL \fall X\in TL.
$$
Setting $\va=-e^{-\tau}$ and $\beta=e^{-\tau}(\eta_h+\bar{t})$, 
we obtain that
$
e^{-\tau}(\mu+\bar{t})=h_*\nabla\va+\beta=\Fes
$
is a Combescure transform of $h$. Moreover, $\mu+t=-\va^{-1}(\Fes+t_T)$.
Thus
$$
f\circ t=h-2\va\|\Fes+t_T\|^{-2}(\Fes +t_T).
$$
Finally, uniqueness of $[(\va,\beta)]$ was observed in Proposition \ref{multi} - $(i)$. \qed\vspace{1.5ex}

\noindent{\bf Remark.} In view of Theorem \ref{thm:main}, the invariance under ${\cal L}$-transformations of the property of admitting a Dupin principal normal with integrable conullity 
can also be derived from the invariance of the $\N$-Ribaucour transformation under $\Lal$-transformations.  The latter is as follows. 
To each $\Lal$-transformation $T$ of $h$ we associate  an $\Lal$-transformation $\tilde T$ for 
$f=\Ral^\N_w(h)$ given~by
\begin{itemize}
\item[(i)] $\tilde T=T$ when $T$ is a conformal map 
of $\R^{n+p}$; 
\item[(ii)] $\tilde T=L_{\hat\xi}$ when $T=L_\xi$ 
for  $\xi\in\N^\perp$, where $\hat\xi=\P_t\xi$. 
\end{itemize}
Then, the correspondence $T\mapsto w^T\in\Sal_\N(T(h))$ 
given in (\ref{ruso}) is such that
\be\label{ruso2}
\tilde T f=\Ral^{\N^T}_{w^T}(Th),
\ee
where $\N^T=T\N$ if $T\in O(n+p)$,
$\N^T=\P^i\N$ in the case of an inversion, and $\N^T=\N$ if $T$ is  either an Euclidean translation, an
homothety or a parallel translation.
\vspace{1ex}

   We conclude this section by showing that the 
$\N$-Ribaucour transformation for holonomic submanifolds 
admits a simple coordinate description in terms of solutions 
of completely integrable first order systems of partial differential equations.  \medskip

Let \mbox{$g\colon L^{n-s}\to\R^{n+p}$} 
be a holonomic submanifold  endowed with principal coordinates 
$(u_1,\ldots,u_{n-s})$, let  $\N$ be a \mbox{parallel} flat normal subbundle, and let $\xi_1,\ldots,\xi_{p+s}$ be a parallel orthonormal normal frame 
such that $\N=\spa\{\xi_{p+1},\ldots\xi_{p+s}\}$. 
Set $ds^2=\sum_{j=1}^{n-s} v_j^2du_j^2,$ 
define $h_{ij}\in C^\infty(M)$ by 
$v_ih_{ij}=\partial v_j/\partial u_i$, $1\le i,j\le n-s$, 
and  $V_j^r\in C^\infty(M)$ by 
\mbox{$v_jA^g_{\xi_r}X_j=V_j^rX_j$}, $1\le r\le s+p$, where $X_j=v_j^{-1}\d/\d u_j$. 
It follows from the Gauss, Codazzi and Ricci equations for $g$  
that the triple $(v,h,V)$ satisfies the completely 
integrable~system of partial differential equations \vspace{1.5ex}
\be\label{fund}
\left\{\begin{array}{l}
\!\!\!(i)\, {\displaystyle\frac{\d v_{i}}{\d u_{j}} = h_{ji}v_{j}},\;\;
(ii)\, {\displaystyle\frac{\d h_{ij}}{\d u_i} + 
\frac{\d h_{ji}}{\d u_j}\! + \!\sum_{k}h_{ki}h_{kj}\! +\! \sum_{r}V_{i}^rV_{j}^r=0},
\vspace{1.5ex}\\
\!\!\!(iii)\; {\displaystyle\frac{\d h_{ik}}{\d u_j} = h_{ij}h_{jk}},\;\;  (iv)\;\displaystyle\frac{\d V_{i}^r}{\d u_j}=h_{ji}V_{j}^r,
\vspace{1.5ex}\\
\end{array}\right.
\ee
where always $i\neq j\neq k\neq i$. Now consider the linear  system of partial differential equations 
of first order
\be\label{lin}
\left\{ \begin{array}{l}
\d\va/\d u_i\; =
v_i\gamma_i\vspace{.5ex}\\
\d\gamma_j/\d u_i =
h_{ji}\gamma_i,\;\;\; i\neq j,\vspace{.5ex}\\
\d \beta_r/\d u_i = - V_i^r\gamma_i.
\end{array} \right.
\ee
System (\ref{lin}) is also completely integrable, the compatibility conditions being satisfied by virtue of (\ref{fund}). 
Moreover,  $(\va,\gamma,\beta)=
(\va,\gamma_1,\dots,\gamma_{n-s},\beta_1\ldots,\beta_{p+s})$ 
is a solution of (\ref{lin}) if and only if the pair $(\va,\beta=\sum_{r=1}^{p+s}\beta_r\xi_r)$, satisfies (\ref{eq:gnorm}). Therefore, $\Ral^\N_{\va,\beta}(g)$ can be parametrized
as $f\colon L^{n-s}\times\R^s\to\R^{n+p}$ given by
\be\label{eq:param3}
f(u,y) = g - 2\va\nu(\sum_{j=1}^{n-s}\gamma_jg_*X_j +
\sum_{r=1}^p\beta_r\xi_r + \sum_{\ell=p+1}^{s+p}(\beta_\ell+y_{\ell-p})\xi_\ell),
\ee
where $\nu^{-1} = \sum_{j}\gamma_j^2+\sum_{r}\beta_r^2+
\sum_{\ell}(\beta_\ell+y_{\ell-p})^2$.\vspace{3ex}

\noindent{\bf\large  Generalized cylinders} 
\vspace{3ex}

    This section is devoted to characterize the class of submanifolds that are $\Lal$-equivalent to the generalized 
cylinders within the class of submanifolds that carry a Dupin principal normal with integrable conullity. Let us first recall their precise definition.\vspace{1ex}

\noindent {\bf Definition.} Let $g\colon\,L^{n-s}\to\Q_\e^{N}$, $\e=0,1,-1$, be an \ii with a  
parallel flat  normal subbundle  $\Nu$ of rank $s$.  The {\it generalized cylinder\/} 
over $g$ determined by $\Nu$ is the $n$-dimensional submanifold parametrized 
by means of the exponential map of  $\Q_\e^{N}$~as
$$
\gamma\in\Nu\mapsto\mbox{exp}^\e_{g(\pi(\gamma))}(\gamma).
$$ 

\vspace{1ex}

   We start by showing that the generalized cylinders are the only submanifolds that carry a relative nullity distribution with integrable conullity.

\begin{proposition}\po\label{nul2}  Let 
$h\colon L^{n-s}\to \Q_\e^{n+p}$ 
be a simply connected submanifold with a parallel 
flat normal subbundle $\N$ of rank~$s$ such that $\E^{\N^\perp}_0=0$. 
Then the generalized cylinder over $h$ determined by $\N$ has relative nullity of constant dimension $s$ and integrable conullity. 

Conversely, any submanifold $f\colon M^n\to \Q_\e^{n+p}$ with 
relative nullity of constant dimension $s$ and integrable conullity 
arises this way locally. That is, $\N=\E_0|_{L}$ is a parallel flat 
normal subbundle of $h=f|_{L}$ for any leaf $L^{n-s}$ of the conullity 
and $f$ is an open neighborhood of $h$ of the generalized cylinder 
over $h$ determined by $\N$.
\end{proposition}

\proof  We argue for the converse in the case where $\e=0$, the proof of the direct statement being straightforward. By Theorem \ref{thm:main}, we have that $\N=\E_0|_{L}$ is a parallel 
flat  subbundle of the normal bundle of 
$h=f|_{L}$ for any given leaf $L^{n-s}$ of the conullity, 
and that $f$ is
an $\N$-Ribaucour transform of $h$ determined by a unique $w\in\Sal_\N(h)$. Since the leaves of relative nullity are 
affine subspaces, then we must have that $\Fes_w\in\N$ by Proposition \ref{multi} - $(v)$, and hence $f$ is  a generalized cylinder in $\R^{n+p}$.
The case $\e=1$ can be easily  reduced to the Euclidean one by 
taking the cone in $\R^{N+1}$ over the submanifold; details are left 
to the reader. The proof of the case $\e=-1$ can be done similarly.\qed\vspace{1ex}

   Now let $f\colon\,M^n\to\R^{n+p}$ be an \ii that carries a Dupin principal normal $\eta$ of multiplicity $s$ and integrable conullity.  By Theorem \ref{thm:main}, there exist a submanifold  $h\colon\, L^{n-s}\to\R^{n+p}$, a parallel flat normal subbundle $\N$ of rank $s$ and an element $w\in\Sal_\N$ such that $f=\Ral^\N_w(h)$. In the following result, we characterize  those $w\in\Sal_\N$ for which  $f$ is  $\Lal$-equivalent to (the stereographic projection of) a generalized cylinder  in $\Q_\e^{n+p}$.

\begin{theorem}\label{2}\po The  following assertions are equivalent: 
\begin{enumerate}
\item[{\em (a)}]  $f=\Ral^\N_w(h)$ for some  ${\cal L}$-trivial (respectively, conformally trivial) $w\in\Sal_\N$;
\item[{\em (b)}] $f$ is ${\cal L}$-equivalent 
{\em (}resp., conformally equivalent\,{\em )}  to 
{\em (}the stereographic projection if $\e\neq 0$ of~{\em )}  
a generalized cylinder in $\Q_\e^{n+p}$.
\end{enumerate}
Moreover, if $h$ is conformally substantial then $\e=\e(w)$ is uniquely determined.
\end{theorem}

\proof  First we give a parametric description of the generalized cylinders in $\Q_\e^N$ as \mbox{$\N$-Ribaucour} transforms.  Consider $\R^N=\{0\}\times\R^N$ inside $\R^{N+1}$ if $\e=1$, or inside the Lorentzian space 
$\Les^{N+1}=\R^{1,N}$ if $\e=-1$, and then take $\Sp^N\subset\R^{N+1}$ and $\Hy^N\subset\Les^{N+1}$. Then, 
the  generalized cylinder in $\Q_\e^N$ over $h$ determined by $\Nu$ can be parametrized as
\be\label{pacy}
\gamma\in\Nu\mapsto h(x) - (1+\e^2) \frac{\e h(x) + \gamma}
{\e+\| \gamma \|^2},\;\;\; x=\pi(\gamma).
\ee

   On the other hand, if $w\in\Sal_\N$ is ${\cal L}$-trivial, it follows from (\ref{fis}) and (\ref{gen})
that $f=\Ral_w^\N(h)$ can be parametrized as $f(u,t)=h_t(u)$, where 
\be\label{uno}
h_t = h - \left(a\|h\|^2 + 2\<h,v\> + c\right)
\frac{ah+v +\delta + t}{\|ah + v +\delta + t\|^2}.
\ee
Moreover, we can assume that $\delta\in\N^\perp$ 
for $\delta_{\N}$ can always be 
canceled by a reparametrization in (\ref{uno}).   
If $h$ is conformally substantial, then $a,v,\delta$ and $c$ 
are  now completely determined up to a common multiplicative constant by 
Proposition \ref{exp}. In particular,
$$
\e(w) = \mbox{sign}\,(ac-\|v\|^2 + \|\delta\|^2)
$$
is well defined.

\begin{lemma}\po\label{epsilon} If $h$ is conformally substantial,
then $\e(w)$ is invariant by conformal transformations and  
parallel translations $L_\xi$ for any $\xi\in {\cal N}^\perp$. 
\end{lemma}

\proof It follows easily from Proposition \ref{trivialt}.
\qed\vspace{1.5ex}

  To prove the equivalence part of the statement, 
we show that the parametrizations (\ref{pacy}) and (\ref{uno}) correspond by 
${\cal L}$-transformations (resp., conformal transformations)  and stereographic projections.
We use Proposition \ref{trivialt} several times without 
further reference. If $a=0$ in (\ref{uno}), we can make 
$a\neq 0$  by a conformal transformation. In fact, 
we may assume 
that  $c\neq 0$. Otherwise, $v\neq 0$ and 
a translation $T_v$ gives $c\neq 0$. Now an 
inversion yields $a\neq 0$.  
By a homothety $H_a$ followed by 
a translation  $T_v$ of (\ref{uno}), we have that
\be\label{choto}
h'_t = h' - \left(\|h'\|^2 + c_1\right)\frac{h' + \delta + t}
{\|h' + \delta + t\|^2},
\ee
where $h'=ah+v$, $h'_t =ah_t+v$ and $c_1=ca-\|v\|^2$.

\begin{claim}\po\label{claim} {\em We may assume 
that $h'+\delta$ 
in (\ref{choto}) is an immersion.
}\end{claim}
\noindent 
The conformal map 
$C=T_{-(c_1+\|q\|^2)^{-1}q}\circ i\circ T_{q}$ for 
$q\in \R^{n+p}$ takes (\ref{choto}) into
$$
\bar h_t=\bar h-(\|\bar h\|^2 + 
c_1(c_1+\|q\|^2)^{-2})
\frac{\bar h+(c_1+\|q\|^2)^{-1}\bar \delta + t}
{\|\bar h+(c_1+\|q\|^2)^{-1}\bar \delta+t\|^2},
$$
where $\bar h=C(h')$ and 
$\bar\delta = \P_C\delta=\delta - 2\<\delta, h'+q\>i(h'+q)$.  
At each point, we obtain using (\ref{sff}) for $P_0=-q$ that 
$A^{\bar h}_{\bar \delta}=
\|h'+q\|^2 A^{h'}_\delta +2\<\delta,h'+q\>I$. Thus,
$$
I-(c_1+\|q\|^2)^{-1}A^{\bar h}_{\bar \delta}=
-(c_1+\|q\|^2)^{-1}\|h'+q\|^2(A^{h'}_\delta -\sigma(q)I), 
$$
where $\sigma(q)= \|h'+q\|^{-2}(c_1-2\<\delta,h'+q\>+\|q\|^2)$. The 
proof of the claim follows from the fact that $\sigma$ is a nonconstant continuous function of $q$. If otherwise,  
$\|h'\|^2 + c_1=0$, 
which is in contradiction with (\ref{choto}) being a parametrization.

A parallel translation $L_{\hat\delta}(f)$ (see (\ref{ruso2})) yields
\be\label{tres}
h_t'' = h'' - \left(\|h''\|^2 + c_2\right)\frac{h'' + t}
{\|h'' + t\|^2},
\ee
where $h''=h'+\delta$, $h''_t=h'_t+\delta_t$ and $c_2=c_1+\|\delta\|^2$.
Now a homothety $H_{|c_2|^{-1/2}}$, if necessary, 
and Lemma \ref{epsilon} yield
\be\label{cuatro}
g_t=g - \left(\|g\|^2 +\e\right)\frac{g+t}{\|g+t\|^2},
\;\;\;\e=\e(w) = 0,\pm 1,
\ee
where $g=|c_2|^{-1/2}h''$ and $g_t=|c_2|^{-1/2}h_t''$ when 
$c_2\neq 0$. 

It remains to show that the stereographic projection 
$$
S=T_{\e e_0} H_{1+\e^2}\, i\, T_{-\e^2e_0}\colon\,\R^N\to\Q^N_\e,
$$ 
where $e_0=(1,0)$, takes the normal form (\ref{cuatro}) 
to a generalized cylinder in $\Q^N_\e$. Notice that  $S=i$ if $\e=0$. 
For each $t\in\Np$, the corresponding parallel
normal vector field $\hat t\in\Nu''$ of 
$k = S(g)=\e e_0 + (1+\e^2)(\| g\|^2 + \e)^{-1}
(g - \e^2 e_0)$ is 
$\hat t = t - 2\langle t, \tilde g\rangle i(\tilde g)$, 
where $\tilde g = g-\e^2 e_0$.
Thus $k_{\hat t} = S(g_t)$ is given by
$$
k_{\hat t}=\e e_0 + (1+\e^2)\frac{\|\tilde g + 
\e^2 e_0 + t\|^{2} \tilde g-\|\tilde g \|^2 (\tilde g + 
\e^2 e_0 +t)} {\e+\| \hat t \|^2}
= k - (1+\e^2) \frac{\e k + \hat t}
{\e+\| \hat t \|^2},
$$
which is, precisely, a generalized cylinder 
in $\Q^N_\e$ over $k$ determined by $\Nu$.\qed\vspace{1ex}

\begin{corollary}\po \label{cor:leaves} Let $f\colon M^n\to \R^{n+p}$ be an \ii 
carrying a Dupin principal normal with integrable conullity.
The following assertions are equivalent:
\begin{itemize}
\item[{\em (a)}] $f$ is $\Lal$-equivalent to (the stereographic projection of)
a generalized cylinder;
\item[{\em (b)}] There is a pair of immersed conullity 
leaves of $f$ that are $\Lal$-equivalent;
\item[{\em (c)}] Any pair of immersed conullity leaves of $f$ 
are $\Lal$-equivalent.
\end{itemize}
\end{corollary}

 The concepts of rotation submanifold and tube admit the 
following extensions for an Euclidean submanifold $g\colon\,S\to\R^N$  with a parallel flat normal 
subbundle $\N$. 
\begin{itemize}
\item[(i)] The {\it generalized rotation submanifold\/} 
$\psi\colon\,\N\to\R^N$ over $g$ determined by $\N$ and $e\in\R^N$ 
is given by
$$
\psi(\gamma)= g(x) - 2\<g(x),e\>\frac{e+\gamma}{\| e+
\gamma\|^2}, \;\;\; x=\pi(\gamma).
$$
\item[(ii)] The  {\it generalized tube\/} $\psi\colon\,\N_1\to\R^N$ over 
$g$ determined by $\N$ and $a\in\R^*$ is given by
$$
\psi(\gamma)=g(x)+a\gamma,\;\;\; 
x=\pi(\gamma),
$$
where $\N_1\subset\N$ denotes the unit sphere subbundle.   
\end{itemize}

\begin{proposition}\po\label{basic} The stereographic 
projection on $\R^N$ of a generalized cylinder in $\Q^N_\e$ 
over $h$ determined by $\Nu$ for $\e\neq 0$ is 
$\Lal$-equivalent 
to one of the following submanifolds:
\begin{itemize}
\item[{\em (a)}] A generalized tube if  $\e=1$ and 
$\Nu$ is not maximal;
\item[{\em (b)}] A generalized rotation submanifold if 
$\e=-1$.
\end{itemize}
Conversely, any generalized tube or rotation submanifold is ${\cal L}$-equivalent to the stereographic projection of a generalized cylinder 
in $\Q^N_\e$ for $\e=1$ or $\e=-1$, respectively.
\end{proposition}

\proof From the proof of Theorem \ref{2}, we know that the stereographic projection of a generalized cylinder in 
$\Q^N_\e$ has the form (\ref{cuatro}), 
with $t\in\tilde\Nu=S(\Nu)$. 
If $\e=-1$, take $e\in\R^N$ such that $\|e\|^2=1$.  
By a translation $T_{-e}$ in (\ref{cuatro}) we have that
$$
g'_t=g'-\left(\|g'\|^2+2\<g',e\>\right)
\frac{g'+e+t}{\|g'+e+t\|^2},
$$
where $g'=g-e$ and $g'_t=g_t-e$. Composing with an 
inversion $i$ yields
$$
i(g'_t)=i(g')-2\left(\<i(g'),e\>+1/2\right)
\frac{e+t}{\|e+t\|^2}.
$$
After a translation $T_{e/2}$, we obtain that  
$h_t=i(g_t')+e/2$ 
is a generalized rotation submanifold over $h=i(g')+e/2$.
\vspace{1.5ex}
  
For $\e=1$ we need an alternative description of 
a generalized tube. 
Take a parallel $\nu\in\N_1$ and set $\tilde\N=\N\cap\{\nu\}^\perp$.
The generalized tube  $f\colon\,\tilde\N\to\R^N$ over $h$ determined by 
$\N$ is given by
$$
f(\gamma') = h - \nu + 2 \frac{\nu+\gamma'}{\|\nu+\gamma'\|^2}.
$$ 
Since $\Nu$ is not maximal, there is a unit parallel 
$\xi'\in \tilde\Nu^\perp$. As in the proof 
of Claim \ref{claim}, we 
use the conformal map $C$ for some $q\in \R^{n+p}$ to replace $\epsilon=1$ in (\ref{cuatro}) by $(1+\|q\|^2)^{-2}$ and to obtain that $C(g)+\xi$ is an 
immersion, where $\xi=(1+\|q\|^2)^{-1}\P_C\xi'$.
Then, the parallel translation $L_{\hat\xi}$ of (\ref{cuatro}) yields
$$
g_t'=g' - \|g'\|^2\frac{g'-\xi + t}{\|g'-\xi + t\|^2},
$$
where $g'=C(g)+\xi$ and $g_t'=C(g_t)+\xi_t$. We obtain a generalized 
tube by composing with the inversion $i$.\qed\vspace{1.5ex}
 
\noindent{\bf Remark.} The three classes in part $(b)$ of Theorem \ref{2} for distinct $\e$ do not have to be disjoint if $h$ 
is not conformally substantial. For example,  tubes may also 
be  rotational submanifolds.  In fact, one can see that an 
element belonging to any two classes also belongs to the 
third.\vspace{1ex}

  As an application of Theorem \ref{2}, we are now able to give a short 
proof of the main result in \cite{dft} with the additional assumption
that the submanifold is locally conformally substantial, thus showing 
the advantage one may have in working with parametric descriptions 
instead of the fundamental equations of the submanifold.

\begin{theorem}\po\label{umb} Let $f\colon\,M^n\to\R^N$ be 
a locally conformally substantial submanifold with a Dupin principal  normal of multiplicity 
$k$ such that its conullity is totally umbilical in $M^n$. If $k=n-1$, 
assume further that the integral curves of the conullity 
are circles in $M^n$. Then $f(M)$ is conformally congruent 
to an open subset of one of the following 
submanifolds:
\begin{itemize}
\item[{\em (a)}] $M^n = L^{n-k}\times\R^k$, and 
$f=(g,id)$ for a submanifold $g\colon\,L^{n-k}\to\R^{N-k}$;
\item[{\em (b)}] $M^n=CL^{n-k}\times\R^{k-1}$, and 
$f=(Cg,id)$ where  $Cg$ is the cone over
a spherical submanifold $g\colon\,L^{n-k}\to\Sp^{N-k}\subset\R^{N-k+1}$;
\item[{\em (c)}] $M^n=L^{n-k}\times_\rho\Sp^k$, and 
$f=(g,\rho\, i)$ for a submanifold 
$g\colon L^{n-k}\to\R^{N-k-1}$, the inclusion $i\colon\,\Sp^k\to\R^{k+1}$ 
and a function $\rho\in C_+^\infty(L)$.
\end{itemize}
\end{theorem}

\proof We use parametrization (\ref{rib2}) and claim that
$w$ is conformally trivial. Since the  nullity distribution 
is $\U(h_t(x))=\P_{[(\va,\beta+t)]}{\N}(x)$, then the 
conullity 
distribution is totally umbilical in the manifold if 
and only if for each $\delta\in {\N}$ 
there exists $\kappa^\delta_t\in\R$ such that
$$
\<A^{h_t}_{\delta_t}X,Y\>=\kappa^\delta_t\<X,Y\>
$$
for all $X,Y\in\U^\perp(t) = T_{h_t}L$. At the leaf
parametrized by $h$ $(\| t \| \to \infty)$, we have that
$A_\delta^h = \kappa^\delta \id$ for all $\delta\in {\N}$.  
We obtain from Theorem \ref{cod} - $(i)$ that
$$
\kappa^\delta_t (\id - 2\varphi\nu_t \Phi_t) 
= \kappa^\delta\id+ 2\<\beta_t,\delta_t\>\nu_t\Phi_t.
$$ 
We easily conclude that $\Phi = a\id$ for some $a\in\R$, 
and the proof 
of the claim follows  using that $\Fes_*=f_*\circ\Phi$. \vspace{.3cm}

\noindent{\it Case $\kappa^\delta= 0$ for all $\delta\in {\N}$}. 
Observe that  the conullity being totally umbilical in the 
manifold is a conformally invariant property. Since ${\N}$ is parallel and totally  geodesic by assumption, we conclude that 
$h$ reduces codimension to $N-k$. Thus, up to 
translation and homothety, we have an orthogonal splitting $\R^{N}=\R^{N-k}\oplus\R^k$ such that $h\subset \R^{N-k}$, ${\N}=\R^k$ and (\ref{cuatro}) takes the form
$$
h_t=(h,0)-\left(\|h\|^2 +\e\right)
\frac{(h,t)}{\|(h,t)\|^2},\;\;\; e = 0,\pm 1.
$$
Choose $e\in\R^N$ such that $e=0, e_N, e_1$ for 
$\e= 0,1,-1$, respectively.  It is easy to see, by composing
with the conformal transformation $T_{-e/2}\,i\,T_e$, 
that we obtain 
cases $(a), (b), (c)$ in the statement for $\e= 0,1,-1$, respectively.
\vspace{1ex}

\noindent{\it Case} $\kappa \neq 0$.
Now $h$ is totally umbilical with respect
to the subbundle $\N$. It is a standard fact (cf.\ \cite{ya}) that $h(L)$ is contained in a sphere $\Sp^{N-k}\subset\R^{N-k+1}$ , which we may assume
to be  centered at the origin,  and that
${\N}=\spa\{\xi\}\oplus {\N}'$, where $\xi$ is the position vector of
$\Sp^{N-k}$ in $\R^{N-k+1}$ and $\N'=\R^{k-1}$  is the orthogonal  complement of $\R^{N-k+1}$ in $\R^N$. 
Now, an inversion with 
respect to a sphere centered at a point in $\Sp^{N-k}$ reduces this 
case to the first one.

Observe now that the three types cannot be glued together by the last
part of Theorem \ref{2}, since $h$ is conformally substantial.
\qed
\vspace{3ex}

\noindent {\bf\large  Weakly reducible Dupin submanifolds}\vspace{3ex}

  Our main goal in this section is to describe how to construct locally all weakly reducible $k$-Dupin submanifolds as defined below. As a consequence, we obtain an explicit coordinate description of a recursive procedure to construct all the holonomic $k$-Dupin submanifolds. Several related results on $k$-Dupin submanifolds are also given.\vspace{1.5ex}

   It is a well-known  fact (see \cite{re1}) that at each point $x\in M^n$ of an \ii  $f\colon\,M^n\to \R^{n+p}$ with flat normal bundle there exist an integer $k(x)$  and unique principal  normals $\eta_1,\ldots,\eta_k\in T_x^\perp M$  such that the tangent space splits orthogonally as 
\be\label{nor2}
T_xM=\E_{\eta_1}(x)\oplus\cdots\oplus\E_{\eta_k}(x).
\ee 
We call $f$ {\it proper\/}  if $k=k(x)$ is constant on $M^n$. 
In this case, each $\eta_j$ is smooth and    
the dimension of $\E_{\eta_j}$ is constant. Hence, $\E^f=(\E_{\eta_1},\ldots, \E_{\eta_k})$ is an orthogonal $k$--net  on $M^n$, that is,  an orthogonal decomposition of $TM$ into $k$ integrable subbundles (cf.\ \cite{rs}). \vspace{1ex}

\noindent{\bf Definition.\/}  An isometric immersion 
$f\colon\,M^n\to\R^{n+p}$ with flat normal bundle is a
\mbox{{\it $k$--Dupin submanifold\/}} if  it is proper and any one of its principal normals $\eta_1,\ldots,\eta_k$ is~Dupin. 
\vspace*{0,1ex}

We now introduce the main concept of this section.
\vspace{1ex}
 
\noindent{\bf Definition.\/} A $k$--Dupin submanifold is {\em weakly reducible\/} if it has a principal normal with integrable conullity. \vspace{1ex}

     Given a $k$--Dupin submanifold $f\colon\,M^n\to \R^{n+p}$, 
we call a Codazzi tensor $\Phi$ on $M^n$ a {\it Dupin tensor adapted to\/} $\E^f$ if  there exist $\phi_1,\ldots, \phi_k\in C^\infty(M)$ such that each function\\ $\phi_j$~is constant along $\E_{\eta_j}$ and $\Phi=\sum_{j=1}^k\phi_j P_{\E_{\eta_j}}$, where $P_{\E_{\eta_j}}$ denotes the orthogonal projection of $TM$ onto $\E_{\eta_j}$.  \vspace{1ex}

\noindent{\bf Definition.\/}  Given a parallel flat normal subbundle $\N$ on a $k$--Dupin submanifold $f\colon\,M^n\to\R^{n+p}$, the \mbox{$\N$-Ribaucour} transform $\Ral^\N_w(f)$ of $f$ 
determined by $w\in\Sal_\N$ is of {\em Dupin type\/} if $\Phi_w$ is a  Dupin tensor adapted to $\E^f$ .\vspace{1ex} 

Finally, given $w=[(\va,\beta)]\in\Sal_\N$ we call $\Ral^\N_w(f)$ {\em regular\/} if $\N\neq 0$  and the vector fields $\bar{\beta}=-\va^{-1}\beta_{\N^\perp}, (\eta_1)_{\N^\perp},\ldots, (\eta_k)_{\N^\perp}$, are 
everywhere distinct. Notice that regularity implies that $\E(w)=0$.\vspace{1ex}

 We are now in a position to prove the main result in this section. The  assumption of regularity in the direct statement 
is only needed to assure that the number of principal normals of the submanifold generated by the \mbox{$\N$-Ribaucour} transformation of a \mbox{$(k\!-\!1)$-Dupin} submanifold is nowhere less than $k$.

\begin{theorem}\po\label{8} Let $h\colon\; L^{n-s}\to\R^{n+p}$ be a simply connected $(k\!-\!1)$--Dupin submanifold with a  parallel flat normal subbundle $\N$ of rank $s$.  Then any regular $\N$-Ribaucour transform of Dupin type  of $h$ is an $n$--dimensional weakly reducible $k$--Dupin submanifold in an open neighborhood of $h$.

Conversely, let $f\colon\,M^n\to\R^{n+p}$ be a weakly 
reducible $k$--Dupin submanifold. Then there exists a \mbox{$(k\!-\!1)$-Dupin} submanifold $h\colon\,L^{n-s}\to\R^{n+p}$ and a parallel flat  normal subbundle $\N$ of $h$ with rank $s$ such that $f$ is locally a regular $\N$-Ribaucour transform of $h$ of Dupin type. \end{theorem}

\proof Let $f=\Ral^\N_w(h)$ be a regular $\N$-Ribaucour transform of Dupin type of $h$. By Proposition \ref{multi} - $(iv)$, the principal normals of $f$ are $\P_t\bar\beta$ and
\be\label{npn}
\bar{\eta}_j
=(\lambda_j^t)^{-1}\P_t(\eta_j
-2\va\nu_t\rho_j^t\bar{\beta})_{\N^\perp},
\ee
where $w=[(\va, \beta)]$, $\eta_1,\ldots, \eta_{k-1}$ are the principal  normals 
of $h$, the $\rho^t_j$  are the eigenvalues of $\Phi_t$ and 
the $\lambda^t_j=1-2\va\nu_t\rho^t_j$ are the ones of $D_t$. Since
$\lim_{t\to\infty}\bar\eta_j =(\eta_j)_{\N^\perp}$,  we 
have from the regularity assumption 
that $\bar{\eta}_1,\ldots,\bar{\eta}_{k-1}$ and $\P_t\bar\beta$
are pairwise distinct on an open neighborhood $U$ of the 
section at infinity of $\N$, and that $\E_{\bar\eta_j}=\E_{\eta_j}$.
A long but straightforward computation shows that
\be\label{ros}
(\lambda_j^t)^2\nap_{X_j}\bar{\eta}_j=
\P_t(\lambda_j^t\nap_{X_j}\eta_j
-2\va\nu_t X_j(\rho_j)(\bar{\beta}-\eta_j)_{\N^\perp}),
\ee
where $X_j\in\E_{\eta_j}$. We have that  $\nap_{X_j}\eta_j=0$ because $h$ is a $(k\!-\!1)$--Dupin submanifold, and that $X_j(\rho_j)=0$ because $\Phi_w$ is a Dupin tensor. It follows from (\ref{ros}) that $\nap_{X_j}\bar{\eta}_j=0$, and hence $\bar{\eta}_j=0$ is a Dupin  principal normal for $1\leq j\leq j-1$.  Moreover, $\P_t\bar\beta$ is also a Dupin principal normal and has integrable conullity by Proposition \ref{1}. Therefore $f|_U$ is an $n$--dimensional weakly reducible \mbox{$k$--Dupin} submanifold.

  Conversely, let $f\colon\,M^n\to \R^{n+p}$ be a $k$--Dupin submanifold and let $\eta_k$ be a principal 
normal of $f$ such that $\E_{\eta_k}^\perp$ is integrable. For 
a leaf $L^{n-s}$ of $\E_{\eta_k}^\perp$, it
is easy to  see that  $h=f|_{{L^{n-s}}}$ has flat normal bundle with principal  normals
$\bar{\eta}_1,\ldots,\bar{\eta}_{k-1}$ given by \mbox{$\bar{\eta}_j=\eta_j+H_j$},
where $H_j$ is the mean curvature vector of $\E_{\eta_j}$. 
Since $\eta_j$
is parallel in the normal connection of $f$, we have that
\be\label{eq:sec}
\nab_{X_j}\bar{\eta}_j
=-\|\eta_j\|^2X_j+\nabla_{X_j}H_j+\alpha_f(X_j,H_j).
\ee
The right hand side of (\ref{eq:sec}) has no
$\E_{\eta_j}^\perp$-component since the last term vanishes and
$\E_{\eta_j}$ is spherical. Hence $h$ is a $(k-1)$--Dupin submanifold. By Theorem~\ref{thm:main}, we have that $\N=\E_{\eta_k}|_{L^{n-s}}$ is a parallel flat normal 
subbundle of $h$ and that $f(M)$ is locally an open neighborhood of a regular $\N$-Ribaucour transform of~$h$. It follows from (\ref{ros}) that $X_j(\rho_j)=0$, $1\leq j\leq k-1$, that is, the
$\N$-Ribaucour transform is of Dupin type. 
\qed\vspace{1.5ex}

 Given a  $k$--Dupin submanifold $f$ and $[(\va,\beta)]\in\Sal_0(f)$, we say that a Ribaucour transform  $\Ral_{[(\va,\beta)]}(f)$ 
is {\em regular\/} if $\lambda_j^{-1}(\eta_j-\bar{\beta})$, $1\leq j\leq k$,  are everywhere nonzero and pairwise 
distinct, where the $\lambda_j$ are the eigenvalues of $D_{[(\varphi,\beta)]}$ and $\bar\beta=-\va^{-1}\beta$. As a consequence of 
the proof of Theorem \ref{8}, we obtain the following characterization of the Ribaucour transformations that 
preserve the class of $k$--Dupin submanifolds.  

\begin{corollary}\po \label{cor:pres} A regular Ribaucour transform of a $k$--Dupin submanifold is also a $k$--Dupin submanifold if and only if it is of Dupin type.
\end{corollary}

\proof We have from (\ref{npn}) that
$\bar{\eta}_i-\bar{\eta}_j=\P(\lambda_i^{-1}(\eta_i-\bar{\beta})
-\lambda_j^{-1}(\eta_j-\bar{\beta})_{\N^\perp})$,
and the result follows.\qed\vspace{1.5ex}

Corollary \ref{cor:pres} generalizes Theorem $2.8$ in \cite{cft}, where it was proved for holonomic Dupin hypersurfaces.  In particular, it shows that the class of $k$-Dupin submanifolds is invariant under $\Lal$-transformations. In view of 
Theorem \ref{2}, we have also the following consequence of Theorem \ref{8}.

\begin{corollary}\po \label{cor:cylinders} A submanifold that 
is ${\cal L}$-equivalent to (the stereographic projection of) 
a generalized cylinder over a submanifold $h\colon\; L^{n-s}\to\Q^{n+p}_\e$ is a $k$--Dupin submanifold if and only if $h$ is a $(k\!-\!1)$--Dupin submanifold and the regularity condition is satisfied.
\end{corollary}

\noindent{\bf Definition.\/} A $k$--Dupin submanifold is
{\em ${\cal L}$-reducible\/} if it \mbox{is 
${\cal L}$-equivalent} to (the stereographic projection 
of) a generalized cylinder over a \mbox{$(k\!-\!1)$--Dupin} submanifold in $\Q^{n+p}_\e$.\vspace{1ex}

By Proposition \ref{basic}, 
the class of ${\cal L}$-reducible $k$--Dupin submanifolds
include the ones that are ${\cal L}$-equivalent to those obtained as in anyone of Pinkall's examples by starting with a 
$(k\!-\!1)$--Dupin submanifold of arbitrary codimension, which we call {\em reducible\/}. Clearly, for Dupin hypersurfaces 
this coincides with the usual notion of reducibility. 
Thus, we have the following 
implications for $k$--Dupin submanifolds; the validity of their converses is discussed 
at the end of this section:
\be\label{chain}\mbox{Reducible} \Longrightarrow 
{\cal L}\mbox{\,-reducible} \Longrightarrow \mbox{Weakly reducible}.
\ee

   One main application of Theorem \ref{8} is for the class 
of holonomic \mbox{$k$--Dupin} submanifolds.  Observe that 
starting  in Theorem \ref{8} with a holonomic $(k\!-\!1)$--Dupin submanifold yields a holonomic $k$--Dupin submanifold, for we have seen that $\E_{\bar\eta_j}=\E_{\eta_j}$.
Conversely, holonomic \mbox{$k$--Dupin} submanifolds are constructed from holonomic $(k\!-\!1)$--Dupin submanifolds.
Therefore, Theorem \ref{8} provides the inductive step for a recursive procedure to construct all holonomic Dupin submanifolds. We  derive next an explicit coordinate 
description of this construction. 

For our construction  we have to use a principal system of coordinates on a holonomic $k$--Dupin submanifold  which
we call a {\it natural coordinate system}. By 
that we mean that the coordinates for each spherical 
leaf of $\E_{\eta_j}$ for $1\leq j\leq k$ are conformal.
In fact, the recursive construction given by Theorem~\ref{8} 
yields such coordinates. To see this,  observe that the parametrization of the leaves of $\E_{\eta_k}$ for the 
generated principal normal $\eta_k$ is conformal by  
\mbox{Proposition \ref{multi} - $(v)$}. Since $D_t|_{\E_{\eta_i}}=\lambda^t_i\id$, $1\le i\le k-1$,  the parametrization of the spherical leaves of  $\E_{\eta_i}$
remains conformal under the transformation.

 Let $h\colon L^{n-s}\to \R^{n+p}$ be a holonomic $(k\!-\!1)$--Dupin submanifold endowed with a natural 
coordinate  system $(u_1,\ldots , u_{n-s})$. 
  For the statements of the next results, we agree that 
$1\le i,j,\ell\le n-s$,  $1\le m\le k-1$ and $1\leq r\leq p$. 
For each index $i$, let $i'$  with $1\le i'\le k-1$ be 
such that 
\mbox{$\d/\d u_i\in \E_{\eta_{i'}}$}. Set \mbox{$v_{i'}=\|\d/\d u_i\|$} and  $h_{jm}=v_{j'}^{-1}\d v_{m}/\d u_j$. 
Given a parallel orthonormal normal frame 
$\xi_1,\ldots, \xi_{s+p}$, we define $V_{i'}^r$ by
$$
A_{\xi_r}\d/\d u_i=v_{i'}^{-1}V_{i'}^r\,\d/\d u_i.
$$  
We call $(v,h,V)$, where $v=(v_1,\ldots ,v_{k-1})$, 
$h=(h_{im})$ and 
$V=(V_m^r)$,  the {\it triple
associated to $h$\/} with respect to the coordinates $(u_1,\ldots , u_{n-s})$ and the normal frame $\xi_1,\ldots, \xi_{s+p}$. We first prove the following fact.

\begin{lemma}\po The triple $(v,h,V)$ satisfies the completely integrable system of partial differential equations
$$
 (I)\;\;\left\{\begin{array}{l}
(i)\; {\displaystyle\frac{\d v_{m}}{\d u_{j}} = h_{jm}v_{j'}},\;\; 
(ii)\; {\displaystyle\frac{\d h_{ij'}}{\d u_i} + \frac{\d
h_{ji'}}{\d u_j} + \sum_{\ell}h_{\ell i'}h_{\ell j'} + \sum_r V_{i'}^rV_{j'}^r=0},
\vspace{1.5ex}\\
(iii)\; {\displaystyle\frac{\d h_{im}}{\d u_j} = h_{ij'}h_{jm}},\;
\;  (iv)\;\displaystyle\frac{\d V_m^r}{\d u_j}=h_{jm}V_{j'}^r, 
\vspace{1.5ex} \\ 
\end{array} \right.
$$
where $\ell'\neq i'\neq j'\neq \ell'$ in $(ii)$ and $i'\neq m$ in $(iii)$.
  Conversely,  let $(v,h,V)$ be a solution of $(I)$ on a simply connected open subset $U\subset\R^{n-s}$  such that $v_m\neq 0$ everywhere. Then there exists a $(k\!-1\!)$--Dupin submanifold $h\colon U\to \R^{n+p}$ such that the 
standard coordinates $(u_1,\ldots , u_{n-s})$ are natural coordinates for 
$h$ and $(v,h,V)$ is the triple associated to $h$ with respect to these coordinates and some parallel orthonormal frame.
\end{lemma}

\proof Equations $(i)$ are merely the definition of $h_{jm}$. From Lemma $1$ 
in \cite{dt1} we have 
\be\label{rib1}
\nabla_{\d/\d u_i}v_{j'}^{-1}\d/\d u_j
=v_{i'}^{-1}h_{ji'}{\d/\d u_i} \fall i\neq j.
\ee
Using this, the remaining equations, except for $(iii)$ when $j'=m\neq i'$ 
and $(iv)$ when $m=j'$, follow by computing the Gauss and Codazzi equations of~$h$. In order to prove that $(iii)$ also holds for $j'=m\neq i'$, 
let $H_m$ be the mean curvature vector of $\E_{\eta_m}$. 
Then we obtain using (\ref{rib1}) that 
\bea
\<\nabla_{\d/\d u_j}H_m,\d/\d u_i\>
\!\!&=&\!\!v_{i'}\d\<H_m,v_{i'}^{-1}\d/\d u_i\>/\d u_j
=-v_{i'}\d (v_m^{-1}h_{im})/\d u_j\\
\!\!&=&\!\! -v_{i'}v_m^{-1}(\d h_{im}/\d u_j-h_{ij'}h_{jm}),
\eea
and our claim follows from the fact that $\E_{\eta_m}$ is spherical. Finally, $(iv)$ for $m=j'$ follows from
$
0=\d(V_{j'}^rv_{j'}^{-1})/\d u_j
=v_{j'}^{-1}\left(\d V_{j'}^r/\d u_j
-V_{j'}^rh_{jj'}\right),
$
where we have used $(i)$ for $m=j'$.

    Conversely, we have from Proposition $3$ in \cite{dt1} that 
there exists a holonomic submanifold $h\colon U\to \R^{n+p}$ such that $(v,h,V)$ is the triple associated to $h$ with 
respect to the standard coordinates 
$(u_1,\ldots , u_{n-s})$  and some parallel orthonormal 
frame. Since $\d(v_{j'}^{-1}V_{j'}^r)/\d u_j=0$ from $(iv)$, 
it follows that $h$ is a $(k\!-1\!)$--Dupin submanifold 
and that the standard coordinates
are natural.\vspace{1ex}\qed

    For our coordinate description  of the holonomic Dupin submanifolds we also need the following fact. 

\begin{lemma}\po\label{le:syst} 
Let $h\colon L^{n-s}\to\R^{n+p}$ 
be a holonomic $(k\!-1\!)$--Dupin submanifold, 
and let $(v,h,V)$ 
be the triple associated to $h$ with respect to natural 
coordinates $(u_1,\ldots , u_{n-s})$ and some parallel orthonormal normal frame. Then the system of total partial differential equations
\be\label{eq:syst}
\d B_m/\d u_j=h_{jm}B_{j'}
\ee
is completely integrable. Moreover, if $(B_1,\ldots ,B_{k-1})$ 
is a solution of {\em (\ref{eq:syst})}, then the system of 
total partial differential equations 
\be\label{eq:system}
\left\{ \begin{array}{l}
 (i)\;\;{\displaystyle\frac{\d\va}{\d u_i} =
v_{i'}\gamma_i},\;\;\;
 (ii)\;\;{\displaystyle\frac{\d\gamma_j}{\d u_i} =
h_{ji'}\gamma_i},\;\;\; i\neq j,\vspace{1.5ex}\\
(iii)\;\; {\displaystyle\frac{\d\gamma_i}{\d u_i}=B_{i'}
-\sum_{j, j'\neq i'}h_{ji'}\gamma_i+\sum_r \beta_r V_{i'}^r}, \;\;\;
(iv)\;\;{\displaystyle\frac{\d \beta_r}{\d u_i} = - V_{i'}^r\gamma_i,}\vspace{1.5ex}\\
\end{array} \right.
\ee
is also completely integrable. 
\end{lemma}

\proof An easy computation shows that the compatibility conditions of (\ref{eq:syst}) follow from $(I)$-$(iii)$. The compatibility conditions of (\ref{eq:system}) can be verified 
by a straightforward computation using 
$(I)$ and (\ref{eq:syst}).\qed\vspace{1.5ex}
  
  In order to simplify the statement of the next result, 
we call a solution $(\va,\gamma,\beta)$ of (\ref{eq:system}) 
{\it generic\/} if the vectors $-\sum_{r=s+1}^{s+p}\va^{-1}\beta_r\xi_r$ and  $\sum_{r=s+1}^{s+p} v_m^{-1}V_m^r\xi_r$, 
$1 \leq m \leq k-1$,   
are~everywhere pairwise distinct.

\begin{theorem}\po  Let $h\colon\, L^{n-s}\to \R^{n+p}$ be a holonomic $(k-1)$--Dupin submanifold endowed with natural coordinates. If $(\va,\gamma,\beta)$ is a generic solution of {\em (\ref{eq:system})},  
then the  map \mbox{$f\colon\,L^{n-s}\times\R^s\to\R^{n+p}$} given by {\em(\ref{eq:param3})} is, at regular points, a holonomic $k$--Dupin submanifold.

 Conversely, if $f\colon\,M^n\to \R^{n+p}$ is a holonomic $k$--Dupin submanifold and $\eta_{\ell}$ is any of its 
 principal normals,  then $h=f|_{L}$ is a holonomic $(k-1)$--Dupin submanifold for any leaf $L^{n-s}$ of its conullity, and there exists a solution $(\va,\gamma,\beta)$ 
of {\em (\ref{eq:system})} such that $f$ can be parametrized 
by {\em(\ref{eq:param3})}.
\end{theorem}

\proof  It is easily seen that $(\va,\gamma,\beta)$ being a solution of (\ref{eq:system}) and $B_1,\!\ldots\!,B_{k-1}$  a solution of (\ref{eq:syst}) 
is equivalent to the tensor $\Phi=\hess \va -A_{\beta}$ being 
a Dupin tensor adapted to $\E_h$. Therefore, $f$ parametrizes the $\N$-Ribaucour transform $\Ral^\N_w(h)$ of Dupin type of 
$h$ determined by $w=[(\va,\beta)]\in\Sal_\N$, where $\beta=\sum_r \beta_r \xi_r$ and 
$\N$ is the parallel flat normal subbundle of $h$ spanned by $\xi_1,\ldots ,\xi_s$. Moreover, the solution $(\va,\gamma,\beta)$ of (\ref{eq:system}) 
being generic is equivalent to the $\N$-Ribaucour transform $\Ral^\N_w(h)$ 
being regular. The result now follows from Theorem \ref{8}.
\vspace{1.5ex}\qed

 In order to derive a sufficient
condition  for a $k$--Dupin submanifold to be holonomic, we define the {\it local conformal codimension} of an \ii $f\colon\,M^n\to\R^{n+p}$ as
$$
c_\ell(f)=\min\{c(f|_U): U\subset M^n \ \mbox{is open}\}.
$$
Recall that $f$ is called {\it $1$-regular} if the  first normal spaces 
$$
N^f_1(x)=\spa\{\a_f(X,Y): X,Y\in T_xM\}
$$
have constant dimension.

\begin{proposition}\po\label{thm:cf} If $f$
is a  $1$-regular $\,$ $k$--Dupin submanifold then
$c(f)\leq k-1$. Moreover, if $c_\ell(f)=k-1$ then $f$ is
holonomic.
\end{proposition}

Proposition \ref{thm:cf} is an easy consequence of the
following  results.

\begin{lemma}\po \label{cor:regu0} Let 
$f\colon\,M^n\to\R^{n+p}$
be a proper \ii with flat normal bundle and principal normals
$\eta_1,\ldots ,\eta_k$. Then $\eta_\ell$ has integrable conullity
if the vectors $\eta_i-\eta_{\ell}$ and $\eta_j-\eta_{\ell}$
are everywhere linearly independent for any pair of indices
\mbox{$1\leq i\neq j\leq k$} with $i,j\neq \ell$.
\end{lemma}

\proof The Codazzi equation implies that
\be \label{eq:codz2}
\<\nabla_{X_i}X_j,X_\ell\> (\eta_j-\eta_{\ell})
=\<\nabla_{X_j}X_i,X_\ell\> (\eta_i-\eta_{\ell})
\ee
for all unit vectors
$X_i\in \E_{\eta_i}$, $X_j\in \E_{\eta_j}$ and
$X_{\ell}\in \E_{\eta_{\ell}}$.\qed

\begin{lemma}\po \label{cor:regu} Let 
$f\colon\,M^n\to\R^{n+p}$
be a proper \ii with flat normal bundle and principal normals
$\eta_1,\ldots ,\eta_k$. At $x\in M^n$ set
$$
S_f(x)=\spa\,\{\eta_i(x)-\eta_j(x) : 1\le i,j\le k\}.
$$
Then $\dim S_f(x)\le k-1$, and  $f$ is holonomic if equality holds
everywhere.
\end{lemma}

\proof The first assertion follows from
$S_f=\spa\{\eta_j-\eta_{\ell},\; 1\le j\le k\}$ for any fixed
$1\leq \ell \leq k$.
If $\dim S_f(x)=k-1$ everywhere, then Lemma \ref{cor:regu0} implies that
the conullity $\E_{\eta_{\ell}}^\perp$ is
integrable for any $1\leq \ell\leq k$, and the second
assertion is a consequence of Theorem~$1$ in \cite{rs}.\qed

\begin{lemma}\po\label{le:sf} Let $f\colon\, M^n\to \R^{n+p}$ 
be a $1$-regular connected $k$--Dupin submanifold with $\dim\,N_1^f=s$. 
Then $f(M)$ is substantially
contained in an affine subspace~$\R^{n+s}$. If $S_f(x)$ has constant dimension,
then  either $S_f=N_1^f$ everywhere or $\dim S_f=s-1$.
Moreover, $\dim S_f=s-1\ge 0$ if and only if $f(M)$ 
is contained in
a sphere $\Sf^{n+s-1}\subset\R^{n+s}$.
In particular, we have  $c(f)=\dim S_f$.
\end{lemma}

\proof  The Codazzi equation  yields 
\be\label{eq:codz}
\nap_{X_j}\eta_i=\<\nabla_{X_i}X_i,X_j\>(\eta_i-\eta_j)
\;\;\;\mbox{if}\;\;\;i\neq j,
\ee
where $X_i\in \E_{\eta_i}$ and $X_j\in \E_{\eta_j}$ are unit
vectors fields. Since $f$ is Dupin,  then the normal vector subbundle $N_1^f$ is parallel
in the normal connection, and the first assertion follows.

 At any point we have that
\be\label{trivia}
\dim N_1^f(x)-1\leq\dim S_f(x)\leq\dim N_1^f(x).
\ee
Assume that $\dim S_f=s-1\ge 0$.  Our claim is trivial for $s=1$,
thus suppose that $s\geq 2$. The principal curvatures corresponding
to a normal vector field $\eta$ are $\<\eta,\eta_j\>$ 
for $1\le j\le k$.
Hence, a smooth unit vector field $\xi$
spanning the orthogonal complement of $S_f$ in $N_1^f$ is an  umbilical
vector field.
For $i\neq j$ we have from (\ref{eq:codz}) that
\bea
0\!\!&=&\!\!\<\nap_{X_j}\eta_i,\xi\>= X_j
\<\eta_i,\xi\>-\<\eta_i,\nap_{X_j}\xi\>
=X_j\<\eta_j,\xi\>-\<\eta_i,\nap_{X_j}\xi\>\\
\!\!&=&\!\! \<\eta_j-\eta_i,\nap_{X_j}\xi\>.
\eea
Thus $\nap_{X_j}\xi\in N_1^f$ is orthogonal to $S_f$, and hence
must vanish.
Therefore, $\xi$ is parallel in the normal connection, and the last assertion follows.\qed \vspace{1.5ex}

\noindent{\it Proof of Proposition \ref{thm:cf}:}
The first claim is an easy consequence of
Lemma~\ref{le:sf} because, if $\dim N_1=k$, then (\ref{trivia}) implies that
$\dim S_f=k-1$ everywhere. By Lemma \ref{le:sf}
the hypothesis on $c_\ell(f)$ now forces $S_f$ to
have constant dimension $k-1$, and the second claim 
follows from Lemma \ref{cor:regu}.
\qed \vspace{1.5ex}

 The next result shows that a $k$--Dupin submanifold must be weakly reducible  if its 
conformal codimension is sufficiently high.

\begin{proposition}\po\label{lll} 
$\!\!$Let $f\colon M^n\to\R^{n+p}$
be a locally weakly irreducible $k$-Dupin \mbox{submanifold}. Then  
$c(f)\le (2/3)k-1$ on each connected component of an open 
dense subset of $M^n$. 
\end{proposition}

\proof  On an open subset $U\subset M^n$ where $f$ is $1$-regular and $S_f$ has constant dimension,   using Lemma \ref{cor:regu0} we have
that for each $1\leq \ell\leq k$ there is an (affine) line 
$L_\ell$ which contains $\eta_\ell$ and at least two more 
principal normals. The estimate for $c(f|_U)$ now follows 
easily from Lemma \ref{le:sf} since $S_f$ is the affine space 
generated by these lines.
\qed\vspace{1ex}

The following example shows that the estimate in the last 
result is sharp.\vspace{1ex} 

\noindent{\bf Example.} 
Take the product immersion of $\ell$ copies of an irreducible isoparametric hypersurface $M^n\subset\Sp^{n+1}\subset\R^{n+2}$ with three distinct principal curvatures. This is a conformally substantial weakly irreducible $n\ell$-dimensional submanifold 
in a sphere of dimension $(n\!+\!2)\,\ell-1$   
for which equality in the estimate  holds.
\vspace{1ex}   

We now discuss whether the converses hold in (\ref{chain}). First we show that the converse is false in the first implication, even for hypersurfaces.

\begin{proposition}\po \label{prop:imp1} There exist $k$-Dupin hypersurfaces with $k\geq 4$ that are ${\cal L}$-reducible but not reducible.
\end{proposition}

\proof Equation (\ref{ruso2}) shows that if a $k$--Dupin hypersurface is weakly reducible with respect to a principal curvature $\lambda$, and $\tilde{\lambda}$ is the corresponding principal curvature of an \mbox{${\cal L}$-transform} of it, then  the conullity leaves of $\lambda$ and $\tilde{\lambda}$ correspond under the \mbox{${\cal L}$-transformation}.  Since the conformal codimension of the conullity leaves of a principal curvature generated by anyone of Pinkall's examples is one and the conformal codimension is invariant under ${\cal L}$-transformations by Proposition \ref{inv}, it follows that all the conullity leaves of a principal curvature of a reducible Dupin hypersurface  have conformal codimension one. Thus, a tube over a weakly irreducible $(k\!-\!1)$--Dupin submanifold $h:L^{n-s}\to\R^{n+1}$ with $k\geq 4$ and $c(h)\geq 2$ is an ${\cal L}$-reducible $k$--Dupin hypersurface that is not reducible. The following well-known fact shows that any irreducible (as a R!
 iemannian manifold) isoparametric submanifold with conformal codimension at least two can be taken as such an $h$.\qed

\begin{proposition}\po \label{prop:isop} Any locally irreducible (as a Riemannian manifold) isoparametric submanifold is weakly irreducible.
\end{proposition}

\proof Let $\eta_1,\ldots,\eta_k$ denote the principal normals of an isoparametric submanifold. 
For any principal normal $\eta_\ell$, the Codazzi equation (\ref{eq:codz}) and the fact that $\eta_\ell$ is parallel in the normal connection imply that $\E_{\eta_\ell}$ is totally geodesic. On the other hand, if $\E_{\eta_\ell}^\perp$ is integrable then  the expressions under parenthesis in (\ref{eq:codz2}) coincide.  Since $k\geq 3$ by the assumption, it follows that both must vanish. Thus $\E_{\eta_\ell}^\perp$ is also totally geodesic, and the de Rham Theorem yields a contradiction. \qed\vspace{1ex}

    For $3$--Dupin hypersurfaces, however,  the three notions of reducibility do coincide. In fact, a weakly reducible 
$3$--Dupin hypersurface can not be Lie equivalent to an isoparametric hypersurface in the sphere by Proposition \ref{prop:isop}, hence the main result of Cecil and Jensen in (\cite{cj}) implies that it must be reducible.

We do not have an explicit example showing that the converse is
false also in the second implication in (\ref{chain}). However, we prove the following result.

\begin{proposition}\po \label{prop:partial} For any $k\geq 4$ there exists a 
holonomic $k$--Dupin submanifold (hence weakly reducible with respect to every principal normal) that is not ${\cal L}$-reducible with respect to some principal normal.
\end{proposition}

By Theorem~\ref{8}, in order to prove Proposition~\ref{prop:partial} it suffices to show that for any $k\geq 4$ there exists a $(k\!-\!1)$--Dupin submanifold $h\colon M^n\to\R^{n+p}$ that carries a nontrivial Dupin tensor adapted to $\E^h$. In fact, we prove that any holonomic $k$--Dupin submanifold $h$ that satisfies $c(h)\leq k-2$ has this property.  This is done by comparing the dimension of the vector space of ${\cal L}$-trivial tensors on $M^n$ with that of Dupin tensors on $M^n$ that are (locally) adapted to $\E^h$. The former is clearly equal to $c(h)+1$ for any $k$--Dupin submanifold $h$.
The latter is computed next for  holonomic submanifolds.

\begin{proposition}\po\label{prop:2} 
Let $f\colon\,M^n\to\R^{n+p}$
be a holonomic $k$--Dupin submanifold. For any $p_0\in M^n$ 
there exist an 
open neighborhood $U$ of $p_0$  and a unique Dupin tensor 
$\Phi$ on $U$ adapted to $\E^f$ such that \mbox{$\Phi(p_0)=\sum_{m=1}^k \phi_m^0 P_{\E_{\eta_m}}(p_0)$} for given
$(\phi^0_1,\ldots, \phi^0_k)\in \R^k$.
In particular, the vector space of Dupin tensors on $U$ adapted to $f$ has dimension~$k$.
\end{proposition}
 
\proof Let $U\subset M^n$  be a simply connected neighborhood 
of $p_0$ endowed with natural coordinates and let 
$\phi_1,\dots,\phi_k$ be smooth functions on $U$. 
It is easily checked that 
the tensor $\Phi=\sum_{m=1}^k \phi_m P_{\E_{\eta_m}}$ is a 
Dupin tensor on $U$ if and only if the functions $B_m=v_{m}\phi_m$ satisfy system (\ref{eq:syst}). The result then follows from the first assertion of Lemma \ref{le:syst}.\vspace{1.5ex}\qed

     We conclude the paper with some consequences of our previous results for $3$--Dupin and $4$--Dupin submanifolds.  

\begin{proposition}\po\label{prop:3dupin} Any nonholonomic $3$--Dupin submanifold is 
Lie equivalent to the stereographic projection of an isoparametric hypersurface in the sphere. 
\end{proposition}

\proof Since any $2$--Dupin submanifold is holonomic, it follows from Theorem \ref{8} that a nonholonomic 
$3$--Dupin submanifold must be weakly irreducible. Moreover, it must also have local conformal codimension one by Proposition \ref{thm:cf}. Therefore, it is (locally) irreducible, and hence Lie equivalent to the stereographic projection of an isoparametric hypersurface in the sphere by the result of Cecil and Jensen \cite{cj}.\vspace{1ex} \qed

     For $4$--Dupin submanifolds the situation is far more complex even globally. As mentioned in the introduction,  there
are examples of compact $4$--Dupin hypersurfaces that are neither weakly reducible nor Lie equivalent to isoparametric hypersurfaces. However, we have the following result for the weakly reducible case.  
 
\begin{theorem}\po Any weakly reducible nonholonomic $4$--Dupin submanifold is ${\cal L}$-equivalent to the stereographic projection of a generalized cylinder over a hypersurface that is Lie equivalent to an isoparametric hypersurface 
in the sphere.
\end{theorem}

\proof Let $h:L^{n-s}\to\R^{n+p}$ be a $3$--Dupin submanifold such that \mbox{$f=\Ral^\N_w(h)$} is not holonomic.  
The Codazzi equation for the Dupin tensor $\Phi_w$ in terms 
of its eigenvalues is 
\be\label{codi}
\left\{\begin{array}{l}
(i)\;\;\, X_j\phi_i+\<\nabla_{X_i}X_j,X_i\>(\phi_i - \phi_j)=0,\vspace{1.5ex}\\
(ii)\;\; \<\nabla_{X_i}X_j,X_k\>(\phi_j - \phi_k) 
= \<\nabla_{X_j}X_i,X_k\>(\phi_i - \phi_k),
\end{array}\right.
\ee
where $X_\ell\in\E_{\eta_\ell}$ and 
$1\le i\neq j\neq k\neq i\le 3$.
Since (\ref{codi}) also holds for any shape operator $A\neq 0$ in the direction of a parallel normal vector field, it follows easily from (\ref{codi})-$(ii)$ and the fact that not all functions $\<\nabla_{X_i}X_j,X_k\>$ can vanish that $\Phi_w=a\id + bA$ for some smooth functions~$a$~and~$b$.  We obtain from (\ref{codi})-$(i)$ that $a,b\in\R$. Therefore
$f$ is ${\cal L}$-equivalent to the stereographic projection of a generalized cylinder over a $3$--Dupin submanifold by Theorem \ref{2}. If such submanifold were holonomic, the same would be true for $f$. The conclusion now follows from Proposition \ref{prop:3dupin}.
\qed

{\renewcommand{\baselinestretch}{1}
\hspace*{-25ex}\begin{tabbing}
\indent \= IMPA -- Estrada Dona Castorina, 110\\  
\> 22460-320 -- Rio de Janeiro -- Brazil\\  
\> email: marcos@impa.br\\
\end{tabbing}}

\vspace*{-7ex}

{\renewcommand{\baselinestretch}{1}
\hspace*{-25ex}\begin{tabbing}
\indent \= IMPA -- Estrada Dona Castorina, 110\\  
\> 22460-320 -- Rio de Janeiro -- Brazil\\  
\> email: luis@impa.br\\
\end{tabbing}}

\vspace*{-7ex}

{\renewcommand{\baselinestretch}{1}
\hspace*{-25ex}\begin{tabbing}
\indent  \=  Universidade Federal de S\~ao Carlos \\
\indent  \= Via Washington Luiz km 235 \\
\> 13565-905 -- S\~ao Carlos -- Brazil \\ 
\> email: tojeiro@dm.ufscar.br \\
\end{tabbing}}

\end{document}